\documentclass[12pt,a4paper]{article}
\usepackage{amsmath,amssymb,latexsym,times,color,fullpage}
\usepackage{helvet,mathrsfs,pifont,bm,tikz}
\usetikzlibrary {shapes}

\def\dd#1{{\,\rm d}#1}
\def\ve{\varepsilon}
\def\rising#1#2{#1^{\overline{{#2}}}}
\def\tr#1{\lfloor#1\rfloor}
\def\qed{{\quad\rule{1mm}{3mm}\,}}
\def\pf{\noindent {\em Proof.}\ }
\def\falling#1#2{#1^{\underline{{#2}}}}
\newcommand{\ds}{\displaystyle}

\newtheorem{thm}{Theorem}
\newtheorem{lmm}{Lemma}
\newtheorem{prop}{Proposition}

\title{\bf Psi-series method in random trees and moments of high
orders}
\author{{\sc Hua-Huai Chern}\\
    Department of Computer Science \\
    National Taiwan Ocean University \\
    Keelung 202 \\
    Taiwan
\and {\sc Hsien-Kuei Hwang}\\
    Institute of Statistical Science\\
    Academia Sinica\\
    Taipei 115\\
    Taiwan
\and {\sc Conrado Mart\'\i nez}\\
    Departament de Llenguatges i Sistemes Inform\`atics \\
    Universitat Polit\`ecnica de Catalunya\\
    Barcelona, E-08034\\
    Spain
}
\date{\today}

\begin{document}
\maketitle

\begin{abstract}
An unusual and surprising expansion of the form
\[
    p_n = \rho^{-n-1}\left(6n +\tfrac{18}5+
    \tfrac{336}{3125} n^{-5}+\tfrac{1008}{3125}\, n^{-6}
    +\text{smaller order terms}\right),
\]
as $n\to\infty$, is derived for the probability $p_n$ that two
randomly chosen binary search trees are identical (in shape and in
labels of all corresponding nodes). A quantity arising in the
analysis of phylogenetic trees is also proved to have a similar
asymptotic expansion. Our method of proof is new in the literature
of discrete probability and analysis of algorithms, and based on the
psi-series expansions for nonlinear differential equations. Such an
approach is very general and applicable to many other problems
involving nonlinear differential equations; many examples are
discussed and several attractive phenomena are discovered.
\end{abstract}



\section{Introduction} \label{sec:intro}
\paragraph{The motivating problem.} This paper was originally
motivated by the following problem. Find the asymptotics of the
sequence $p_n$ defined recursively by
\begin{align}
    p_n = n^{-2}\sum_{0\le j<n}p_jp_{n-1-j}\qquad(n\ge1).
    \label{pn}
\end{align}
with the initial condition $p_0=1$. The sequence $p_n$ is nothing
but the probability that two randomly chosen binary search trees
(BSTs) of size $n$ are identical (having exactly the same shape and
with the same labels for corresponding nodes), and was first studied
by Mart\'\i nez in \cite{Martinez91} as an auxiliary function for
understanding the typical performance of the equality test of two
random BSTs; see below for more background details. A minor
variation of this sequence was encountered in the analysis of
maximum agreement subtrees in \cite{BMS03} under the Yule-Harding
model.

While shape parameters defined on
a single random tree has been extensively studied in the literature
for many varieties of trees, properties of statistics defined on a
pair or multiple of random trees received comparatively less
attention, partly because of the intrinsic complexity of the
underlying analytic problems. Yet many practical situations (such as
\emph{tanglegrams}) naturally lead to such a study, typical example
being the so-called ``hereditary properties" or ``recurrent
properties", which in turn cover the equality, root occurrence,
simplification rules, reduction rules, ``clashes" as special cases;
see \cite{Martinez91,SF06,FS07} for more details.

Recently, there has been more study of statistics defined on two
random combinatorial objects; see \cite{BF09} and the references
therein.

\begin{figure}[h!]
\begin{center}
\begin{tikzpicture}
[level distance=10mm, every node/.style={thick, color = black, text
centered, fill=red!50,ball color=gray!50, circle,
nodes={fill=gray!50,draw=blue,ball
color=gray!50,circle}}, level 1/.style={thick,sibling
distance=30mm}, level 2/.style={thick,sibling distance=15mm},
scale=1.6,minimum height=.5cm,scale=.6]
{\tiny \node {\textbf{6}}
    child {node {\textbf{2}}
        child {node {\textbf{1}}}
        child {node {\textbf{4}}
            child {node {\textbf{3}}}
            child {node {\textbf{5}}}
        }
    }
    child {node {\textbf{8}}
        child {node {\textbf{7}}}
        child {node {\textbf{10}}
            child {node {\textbf{9}}}
            child[missing]
        }
    };}
\end{tikzpicture}\qquad
\begin{tikzpicture}[node distance=60pt]
\node at (4,0) {{\small$\displaystyle\mathbb{P}(U_n=j)=\frac1n$}};
\node at (4,-.8) {{\small$j = 0, \cdots , n-1$}};%
\tikzstyle{every node}=[thick, color = black, text centered ] \node
(A) [fill=gray!50,draw=blue!30!black,ball color=gray!50,circle] {};
\node (B) [below of=A,xshift = -25pt,isosceles triangle,shape border
rotate=90,fill=gray!60,minimum height=2cm,draw=blue!30!black]
{\tiny{$U_n$}};%
\node (C) [right of=B,yshift = -8pt,isosceles triangle,shape border
rotate=90,fill=gray!40,minimum height=2cm,draw=blue!30!black]
{\tiny{$n-1-U_n$}};%
\path [thick,shorten >=1pt] (A) edge (B.north);%
\path [thick,shorten >=1pt] (A) edge (C.north);
\end{tikzpicture}
\end{center}
\caption{\emph{Left: the BST constructed from the sequence
$\{6,2,4,8,7,1,5,3,10,9\}$. Right: the root assumes the value $j+1$
with equal probability $1/n$ for $j=0,\dots,n-1$.}} \label{fig0}
\end{figure}
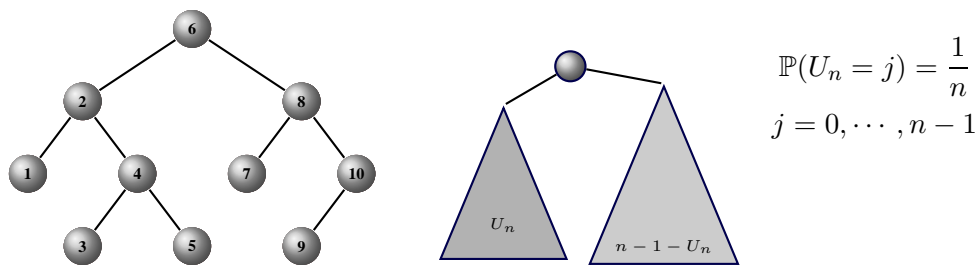

\paragraph{Random BSTs.} For completeness, we first describe BSTs.
Given a sequence of distinct numbers $\{x_1, \dots, x_n\}$, we can
construct the corresponding BST as follows. If $n=0$, then the tree
is empty. If $n\ge1$, then we place $x_1$ at the root; the remaining
numbers are compared one after another with $x_1$, and are directed
to the left subtree of the root if they are smaller, to the right
subtree if larger. Numbers directed to each subtree are constructed
recursively by the same procedure according to their original order;
see Figure~\ref{fig0} for a plot.

By \emph{random BSTs}, we assume that all $n!$ permutations of $n$
distinct elements are equally likely, and construct the BST from a
random permutation. Then we see that the root assumes the value $j$
with probability $1/n$ for $j=1,\dots, n$, which is also the
probability that the left subtree of the root has size $j-1$.

\medskip

\noindent\textbf{Definition: [Equality of two ordered, labeled
trees].} Two ordered, labeled trees of the same size total number of
nodes) are said to be \emph{equal} or \emph{identical} if either
both trees are empty or they have common root label with all
corresponding ordered subtrees equal.

\medskip

The definition extends to the equality of $d$ trees with $d\ge2$.

Now we take two random BSTs independently, and our $p_n$ gives the
probability that the two trees are identical. Equivalently, we take
two random permutations of $n$ elements; then $p_n$ denotes the
probability that the BSTs constructed from these two permutations
are equal. (A simple example: $(2,1,3)$ and $(2,3,1)$ lead to the 
same BST of the shape \raisebox{-.3cm}{\begin{tikzpicture}[scale=1]
\draw[line width=0.8pt] (.38,-.38) -- (.13,-.13);%
\draw[line width=0.8pt] (-.38,-.38) -- (-.13,-.13);%
\draw (0,0) circle (5pt);%
\draw (-0.5,-0.5) circle (5pt);%
\draw (0.5,-0.5) circle (5pt);%
\node  at (0,0) {\tiny{$\mathbf2$}};%
\node  at (-0.5,-0.5) {\tiny{$\mathbf1$}};%
\node  at (0.5,-0.5) {\tiny{$\mathbf3$}};%
\end{tikzpicture}}.)

\begin{figure}
\begin{center}
\begin{tikzpicture}[xscale=.05,yscale=1]
\foreach \x/\xtext in { 10 , 30, ..., 90} \draw (\x,1pt) -- (\x
,-1pt) node[anchor=north] {$\xtext0$}; \foreach \y/\ytext in { 1 ,
2,3} \draw (1,\y ) -- (-1,\y ) node[anchor=east] {$\ytext$};
\draw[-,line width=.5pt] (-1,0) -- (100,0) node[right] {};%
\draw[-,line width=.5pt] (0,-.05) -- (0,3.5) node[above] {};%
\draw[line width=.7pt,smooth,color=red] plot
coordinates{(0.1,1)(0.2,1.414)(0.3,1.651)(0.4,1.824)(0.5 ,
1.955)(0.6,2.059)(0.7,2.143)(0.8,2.214)(0.9,2.274)(1,2.325)
(1.1,2.37)(1.2,2.41)(1.3,2.445)(1.4,2.476)(1.5,2.505)(1.6,2.531)
(1.7,2.554)(1.8,2.576)(1.9,2.596)(2
, 2.614)(2.5,2.688)(3 ,2.743)(4,2.817)(5,2.867)(6,2.902)(7,2.928)(8,
2.949)(9,2.966)(10 ,2.98)(11,2.992)(15
,3.025)(20,3.049)(25,3.064)(30 ,3.075)(35,3.083)(40
,3.089)(45,3.094)(50,3.098)(55 ,3.102)(60 ,3.104)(65
,3.107)(70,3.109)(75 ,3.111)(80 ,3.112)(85
,3.114)(90,3.115)(95,3.116)};%
\end{tikzpicture}
\begin{tikzpicture}[xscale=0.025,yscale=55]
\foreach \x/\xtext in { 50 ,100, ..., 200} \draw (\x ,3.091) -- (\x
,3.089) node[anchor=north] {$\xtext$}; \foreach \y/\ytext in {3.10,
3.12, 3.14} \draw (0.1,\y ) -- (-.1,\y ) node[anchor=east]
{$\ytext$};
\draw[-,line width=.5pt] (-1,3.09) -- (210,3.09) node[right] {};%
\draw[-,line width=.5pt] (0,3.089) -- (0,3.155) node[above] {};%
\draw[line width=.7pt,color=red] plot coordinates{(1, 3.09838)(2,
3.14812)
(3, 3.13990)(4, 3.14098)(5, 3.14083)(6, 3.14086)(200, 3.14086)};%
\end{tikzpicture}
\end{center}
\caption{The figures of $-(\log p_n)/n$ (left) and
$-\log(p_n/(6n+18/5))/(n+1)$ (right).} \label{fg-pn}
\end{figure}
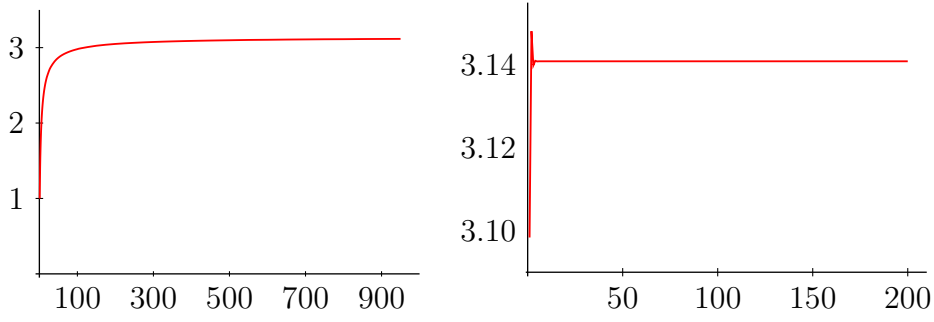

\paragraph{A simple upper bound.} The simple-looking recurrence
(\ref{pn}) can be quickly estimated by the following inductive
argument. If we assume the form $p_n\le c(n+1)\varrho^{-n-1}$ for
$n\ge0$, then we see by induction that
\begin{align*}
    p_n \le \frac{c^2}{n^2}\varrho^{-n-1} \sum_{0\le j<n}
    (j+1)(n-j) = \frac{c^2(n+1)(n+2)}{6n}\varrho^{-n-1}.
\end{align*}
In order that the rightmost term is less than $c(n+1)
\varrho^{-n-1}$, we can take a positive integer $n_0$, let $c :=
6n_0/(n_0+2)$, and then choose $\varrho$ as
\[
    \varrho := \min_{0\le j\le n_0} \left(\frac{6n_0(j+1)}
    {p_j(n_0+2)}\right)^{1/(j+1)}.
\]
Then we obtain
\begin{align} \label{pn-induction}
    p_n \le \frac{6n_0}{n_0+2} (n+1)\varrho^{-n-1} ,
\end{align}
for all $n\ge0$. This gives successively improving bounds for
$\varrho$ for increasing values of $n_0$; see Table~\ref{tab1},
where we take only the first four digits after the decimal point
without rounding. In particular, taking $n_0=6$ leads to the bound
$p_n \le \frac32(n+1) 3^{-n}$.
\begin{table}[h!]
\begin{center}
\begin{tabular}{|c|c|c|c|c|c|c|c|c|c|}\hline
$n_0$& $1$ & $2$ &$3$ &$4$ &$5$ &$6$ &$7$ &$8$ &$9$ \\ \hline%
$\varrho$&$2$ & $2.4494$ & $2.6832$ & $2.8284$ & $2.9277$ & $3$ &
$3.0274$ & $3.0488$ & $3.0659$ \\ \hline\hline $n_0$ & $10$ &$20$
&$30$ &$40$ &$50$ &$60$ &$70$ &$80$ &$90$ \\ \hline
$\varrho$&$3.0794$ & $3.1235$ & $3.1328$ & $3.1362$& $3.1378$ &
$3.1387$ & $3.1393$ & $3.1396$ & $3.1399$ \\\hline
\end{tabular}
\end{center}\caption{\emph{Numerical values of} $\varrho$.}
\label{tab1}
\end{table}
The simple bound (\ref{pn-induction}) obtained by induction and
numerical evidence suggest the possibility that $p_n \sim 6n
\rho^{-n-1}$ for some values of $\rho\approx 3.14$ (see
Figure~\ref{fg-pn}). How to prove this? And is $\rho=\pi$?

\paragraph{The nonlinear differential equation.} As the elementary
argument we used above is not strong enough to derive more precise
asymptotic approximations to $p_n$, we consider instead the
generating function $P(z) := \sum_{n\ge0} p_n z^n$, which satisfies
the nonlinear differential equation (abbreviated throughout as DE)
\begin{align} \label{bst-de}
    zP''(z) + P'(z) = P^2(z),
\end{align}
with the initial conditions $P(0)=P'(0)=1$. This nonlinear DE is of
Emden-Fowler type for which there is no explicit closed form
solution; see \cite{PZ95}. In addition to the apparent singularity
determined by the equation, the DE (\ref{bst-de}) also has
singularities determined by the initial conditions, which are often
referred to as the \emph{movable singularities}.

\paragraph{Frobenius method.} Starting from the DE (\ref{bst-de}),
the next step is often to apply the Frobenius method (see
\cite{Ince26}), namely, we assume the solution of $P(z)$ to be of
the form
\begin{align}
    P(z) = \sum_{j\ge0} c_j (1-z/\rho)^{j-\alpha},
    \label{frobenius}
\end{align}
for some $\alpha$ and $\rho>0$, substitute this form into
(\ref{bst-de}), and then determine $\alpha$ and the coefficients
$c_j$ inductively one after another. This classical procedure yields
$\alpha=2$, $c_0 = 6/\rho$,
\begin{align} \label{c1-c5}
    c_{1}=-{\frac {12}{5\rho}}, \; c_{2}=-{\frac
    {7}{25\rho}},\; c_{3}=-{\frac
    {14}{125\rho}},\; c_{4}=-{\frac {63}{1250\rho}}, \;c_{5}
    =-{\frac {161}{9375\rho}}.
\end{align}
But then \emph{inconsistency} arises since the coefficient of
$(1-z/\rho)^4$ on
\begin{align} \label{c6}
    \text{LHS of (\ref{bst-de})} = \rho^2\left(12c_6 +
    \frac{483}{3125}\right)\not=
    \text{RHS of (\ref{bst-de})} = \rho^2\left(12c_6 +
    \frac{77}{625}\right),
\end{align}
and $c_6$ cannot be determined by simply matching the coefficients
of both sides. This trial suggests that the local expansion of $P$
near the singularity $\rho$ will not be of the form
(\ref{frobenius}) and means that the classical Frobenius method
fails for the nonlinear DE (\ref{bst-de}).

\paragraph{Psi-series method.} We will introduce a different
type of expansion called \emph{psi-series expansion} (or Painlev\'e\
expansion; see \cite{Hille76}) and it will turn out that $P(z)$
admits an asymptotic expansion of the form
\begin{align} \label{UZ}
    U(Z) := \sum_{j\ge0} Z^{j-2}\sum_{0\le \ell \le
    \tr{j/6}} c_{j,\ell}(\log Z)^\ell, \qquad Z := 1-z/\rho,
\end{align}
when $z$ lies near the singularity $\rho$. This form, first
conjectured by Mart\'\i nez in \cite[Ch.\ 9]{Martinez92}, also
explains why the expansion (\ref{frobenius}) leads to inconsistency.
Thus $z=\rho$ is not a pole but instead \emph{a pseudo-pole}; see
\cite{Hille76}. The first few terms of $U(Z)$ are given as follows.
\begin{align} \label{UZ-exp}
\begin{split}
    \rho U(Z) &=6\,Z^{-2}-{\frac {12}{5}}\,Z^{-1}
    -{\frac {7}{25}}-{\frac {14}{125}}\, Z
    -{\frac {63}{1250}}\,Z^{2}
    -{\frac {161}{9375}}\,Z^{3}\\
    &\quad +\rho c_6Z^{4}+\rho\sum_{j\ge 7}
    \sum_{0\le \ell \le\lfloor j/6\rfloor}
    c_{j,\ell} Z^{j-2}\log^{\ell} Z,
\end{split}
\end{align}
for $Z$ small, where $c_6 := c_{6,0}$ and the $c_{j,\ell}$'s are
polynomials of the parameter $c_6\rho$ with degree $\lfloor
(j-6\ell)/6\rfloor$ for $j\ge7$.

The approach we use in this paper is roughly as follows. After
checking the failure of Frobenius method, we construct a suitable
psi-series $U(Z)$ (by matching coefficients) so that $U$ satisfies
\emph{formally} the DE (\ref{bst-de}). The series in (\ref{UZ}) is
\emph{a priori} an asymptotic expansion, but we will show that it is
indeed absolutely convergent in the cut-disk $|Z|\le 1-\ve$,
$Z\not\in[-1+\ve,0]$. Thus the function $U$ is well defined there
and satisfies the DE (\ref{bst-de}) and differs from $P$ only by
their initial conditions. Such a procedure still leaves undetermined
two important parameters (similar to the initial conditions of the
DE (\ref{bst-de})), one is obviously $\rho$ and the other implicit
one is $c_6 := c_{6,0}$ due to the same reason as the Frobenius
method. This means that $U$ is not only
a function of $Z$, but also a function of $\rho$ and $c_6$.

Now to fix $U$ in a unique way, we connect $P(z)$ and $U(Z)$ by
first choosing a number $z_0\in[\ve\rho,\rho-\ve]$, and by
considering the solution $(\rho,c_6)$ of the two equations
\begin{align} \label{sys-eq}
    \left\{\begin{array}{l}
        U(Z_0)= P(z_0) \\
        U'(Z_0)= -\rho P'(z_0),
    \end{array}\right.
\end{align}
where $Z_0 := 1-z_0/\rho$. We will show below
(Proposition~\ref{prop-ac}) that, as a function of $Z$ (or $\rho$)
and $c_6$, the series $U$ has a nonzero radius of convergence for
each finite $c_6$. Also we can easily derive simple upper and lower
bounds for $\rho$ as above. Thus, as a standard initial-value
problem, the system of equations (\ref{sys-eq}) has a unique
solution pair of $(\rho,c_6)$. This determines uniquely the pair
$(\rho, c_6)$. Furthermore, $P$ and $U$ have a common region of
analyticity, and we see by analytic continuation that $U$ is the
exact and asymptotic solution we have been looking for.

Although no analytic forms for $\rho$ and $c_6$ are available, we
can compute the numerical values of $\rho$ and $c_6$ as follows.
First, the values of $U(Z_0)$ and $U'(Z_0)$ can be well approximated
by their partial sums since the terms of the series converge in an
exponential rate; see (\ref{fk-est}); similarly, the values of
$P(z_0)$ and $P'(z_0)$ can be computed by first computing $p_n$ by
its defining recurrence and then summing a sufficiently large number
of initial terms up, the convergence rate being also exponential.
Then we solve successively the corresponding system of equations by
using an increasing number of terms in the partial sums; see next
section for details.

\paragraph{Asymptotics of $p_n$.} From the expansion (\ref{UZ}) and
suitable analytic continuation to be clarified below, we deduce our
main result for $p_n$.
\begin{thm} The probability $p_n$ that two randomly chosen binary
search trees of $n$ nodes are equal satisfies the asymptotic
expansion
\begin{align} \label{pnn}
    p_n \sim \rho^{-n-1}\left(6n +\frac{18}{5} + \sum_{j\ge6}
    n^{-j+1} \sum_{0\le \ell<\tr{j/6}} C_{j,\ell} (\log n)^\ell\right),
\end{align}
for explicitly computable constants $C_{j,\ell}$, where $\rho=
3.14085\,75672\,02936\,95160\dots$
\end{thm}
Thus $\rho\not=\pi$. In particular, the first few terms read
\begin{align*}
    p_n &=\rho^{-n-1}\Biggl(6n +\frac{18}{5}+\frac{336}{3125\, n^5} +
    \frac{1008}{3125\, n^6}+\frac{10416}{15625\, n^7}\\
    &\qquad +\frac{91728}{78125\,n^8}+\frac{8234352}{4296875\,n^9}
    +\frac{12228048}{4296875\,n^{10}}\\
    &\qquad+\frac1{n^{11}}\left(\frac{9483264}{5078125}\,H_n
    +\frac{5621191632}{726171875}+\frac {677376}{1625}c_6\right)
    +O\left(\frac{\log n }{n^{12}}\right)\Biggr),
\end{align*}
where $H_n := \sum_{1\le j\le n}j^{-1}$, and we see that no terms of
the form $cn^{-j}$ with $j=1,\dots,4$ appear in the expansion.
Numerically, the parameter $c_6$ can be determined approximately as
$c_6=- 0.00150\,84982\,09405\,93425\dots$; see the numerical 
discussions on Page~\pageref{page-na} for details.

As far as we were aware, the asymptotic expansion (\ref{pnn}) with
missing terms is rare in the analysis of algorithms and applied
probability literature. The expansion also indicates that the
approximation of $p_n\rho^{n+1}$ by the first two terms $6n+18/5$ is
numerically very precise as can be seen in Figure~\ref{fg-pn}.

\paragraph{Features.} In addition to the unusual form of (\ref{pnn})
and its theoretical value \emph{per se}, the interest of such a
psi-series expansion is multifold. First, since no analytic form for
the movable singularity $\rho$ is available, the psi-series
expansion provides an effective means for obtaining an approximate
value to $\rho$ by the argument we mentioned above; see (\ref{z0z1})
below for more numerical details. Second, from a methodological
point of view, the method of proof we use to prove (\ref{pnn}) is of
some generality. Note that the first two terms on the right-hand
side of (\ref{pnn}) can be easily obtained by the method of matched
coefficients once we assume that $p_n$ has the form (\ref{pnn}).
Third, the precise approximation we derive has direct consequences
in the original motivating problem, as well as several others in the
examples we discuss below. Fourth, such a consideration leads to
several interesting and unexpected phenomena as we will see in 
the following sections.

\paragraph{Outline of this paper.} We describe the psi-series method
and give the proof of the asymptotic expansion (\ref{pnn}) in the
next section. Then we extend in Section~\ref{sec:rts} the
consideration of the probability of equality to either more than two
random BSTs or to other variants of BSTs. It turns out that the
forms of the asymptotic expansion for the probability of equality of
$d$ random BSTs differ drastically according to  the parity of $d$,
a result not intuitively obvious. Section~\ref{sec:2mary} considers
the case of two random $m$-ary search trees and we will see that the
number of missing terms in the asymptotic expansion increases as $m$
grows. Equality of two random fringe-balanced BSTs is considered in
Section~\ref{sec:2fbbst} and there, unlike $m$-ary search trees, the
error term beyond the constant term in the asymptotic expansion does
not change with the structural parameter once it exceeds one,
another unexpected result. Asymptotics of higher-order moments will
then be considered in Section~\ref{sec:mho} with a few
representative examples taken from the cost of partial-match queries
in random trees, random partition structures and solutions of
Boltzmann equations (from statistical physics). We group the details
of some proofs in Appendix.

\paragraph{Notations.} For each problem studied, $\rho$ always denotes
the dominant singularity of the associated nonlinear DE and $Z :=
1-z/\rho$. The symbols $c, c', c_j, c_j', c_{i,j}, C, C_j, C_j',
C_{i,j}, K, K'$ all denote suitably chosen constants, not
necessarily the same at each occurrence.

\section{Psi-series method}
\label{sec:bst}

We discuss in details the psi-series solution to our nonlinear DE
(\ref{bst-de}) and the tools needed to justify it, then we prove
(\ref{pnn}).

\paragraph{Analytic properties of $P(z)$.} First, the solution $P(z)$
to the DE (\ref{bst-de}) has positive radius of convergence and is
analytic at the apparent fixed singularity $z=0$ by definition. By
simple induction as we discussed in the introduction 
(Section~\ref{sec:intro}) and  Pringsheim's
theorem (since all coefficients $p_n$ are positive; see \cite[p.\
240]{FS09}), we expect that $P(z)$ has a finite movable singularity
at, say $z=\rho$, and the asymptotics of $p_n$ will be dictated by
the local asymptotic expansion of $P(z)$ as $z\sim \rho$.

Mart\'{i}nez \cite[p.\ 117]{Martinez92} proved that the function
$P(z)$, originally defined only inside the disk $|z|<\rho$ can be
analytically continued to the cut-disk $|z|\le \rho+\ve\setminus
[\rho,\rho+\ve]$ with $\rho$ being the sole singularity there.

From a theoretic point of view, the movable singularity $\rho$ for
the DE (\ref{bst-de}) can be either of the following types:
\begin{itemize}
\item poles,
\item branch points (algebraic or logarithmic),
\item essential singularity.
\end{itemize}
Simple poles and algebraic points are first excluded because of the
above trial via Frobenius method. We then show that $P$ can be
analytically continued into a function defined by a series expansion
of the form (\ref{UZ}) that converges absolutely in the cut-region
\begin{align} \label{CR}
    \mathscr{C}_R:=
    \{z\,:\,0<|z-\rho|\le R,z\not\in[\rho,\rho+R]\},
\end{align}
for some $R>0$. Thus the possibility that $\rho$ is an essential
singularity is further excluded, and $\rho$ is a logarithmic branch
point (or called pseudo-pole).

Our first focus in this paper is on the determination of the right
form of the solution to (\ref{bst-de}). More detailed and complete
introduction and discussions on the theory related to {\em
Painlev\'{e} analysis} can be found in \cite{CM08,CR00} and
the references therein.

\paragraph{The ARS method (Type checking).} A widely used procedure
to check the singularity type (and the local expansion) of nonlinear
differential equations is the following procedure, often called the
ARS algorithm due to Ablowitz, Ramani and Segur \cite{ARS80}, which
bears some resemblance to the Frobenius method.

In this method, we start assuming that the solution to the DE 
(\ref{bst-de}) admits the formal Laurent expansion (\ref{frobenius}) 
about the cut-disk $\mathscr{C}_R$ for some positive number $R$.
\begin{itemize}

\item[\ding{182}] \emph{Leading order analysis}: Assume $P(z)\sim
c_0(1-z/\rho)^{-\alpha}$. By balancing the dominant terms $\rho
P''(z)$ and $P(z)^2$ in (\ref{bst-de}), we see, as in Frobenius
method, that $\alpha=2$ and the companion constant $c_0=6/\rho$.
Thus we can exclude the possibility of an algebraic singularity.

\item[\ding{183}] \emph{Resonance analysis}: Starting from this pair
$(\alpha,c_0)=(2,6/\rho)$, if the solution admits only poles, then
by substituting (\ref{frobenius}) into (\ref{bst-de}) and by
equating coefficients, the coefficients $c_j$'s are characterized by
the recurrence relation of the form
\begin{align} \label{frob-coeff}
    \Phi (j) c_j =(j-3)^2c_{j-1}+\rho\sum_{1\le j<n} c_j c_{n-j}
    =:G_j (\rho, c_0, c_1, \ldots , c_{j-1} ), \qquad j \ge 1 ,
\end{align}
where $\Phi (j)= (j+1)(j-6)$ and $c_j=0$ for all $j<0$. The roots of
$\Phi (j)$ are called {\em resonance} and $-1$ is always a root of
$\Phi(j)$, reflecting the \emph{arbitrariness of the movable
singularity} $\rho$. For most of our purposes, a less involved and
very commonly used technique is to substitute the test function
\begin{align*}
    c_0(1-z/\rho)^{-\alpha}+c_r(1-z/\rho)^{r-\alpha}
\end{align*}
into the DE (\ref{bst-de}) instead. By collecting the coefficients
corresponding to the term $c_r(1-z/\rho)^{r-4}$, we still get the
same $\alpha, c_0$ and $\Phi(r)$. In this case, we see that $\Phi$
has only one positive resonance $6$ that needs to be further
examined.

\item[\ding{184}] \emph{Compatibility}: Once we have the system
(\ref{frob-coeff}) and identify the resonance, the next step is to
consider its solvability. Obviously, (\ref{frobenius}) is the
solution to (\ref{bst-de}) if and only if all the coefficients
$c_k$'s can be computed recursively by (\ref{frob-coeff}). This fact
defines the \emph{compatibility of the resonance}: for any resonance
$r$ of $\Phi$, if $G_r (\rho, c_0, c_1$, $ \ldots, c_{r-1} ) =0 $ is
satisfied, then the resonance $r$ is said to be compatible;
otherwise, $r$ is incompatible.

From (\ref{c1-c5}) and (\ref{c6}) it follows that $r=6$ is
incompatible. The formal series solution by introducing suitable
logarithmic terms starting at the index $6$ has to be considered
instead (see (\ref{UZ-exp})). The movable singularity $\rho$ to
(\ref{bst-de}) is proved to be a logarithmic branch point since we
will show that the associated series solution is absolutely
convergent in the region $\mathscr{C}_R$ for some $R>0$.
\end{itemize}

In cases when the compatibility of resonance is consistent, the
solution of Laurent expansion is the one we need if it has a
positive radius of convergence. The above ARS Algorithm is useful in
determining if a nonlinear ODE admits the \emph{Painlev\'{e}
property}, namely, the DE has only solutions free from movable
branch points. In our case, the DE (\ref{bst-de}) does not satisfy
the Painlev\'e\ property.

\paragraph{Our approach vs the ARS algorithm.} The method of proof
we use does not, however, rely completely on this method for two
reasons. First, it requires the \emph{a priori} information that
$\rho$ is not an essential singularity, a property often hard to
prove. Second, even we can prove that the singularity is not
essential, the incompatibility of a resonance (or several) may in
some cases very difficult to establish due to the variation of an
additional parameter as in the cases of $d$
random BSTs (Subsection~\ref{sec:d-bst}) and $m$-ary search trees 
(Subsection~\ref{sec:2mary}).

On the other hand, the ARS algorithm does provide an effective means
of computing the exact form of the psi-series expansion for all the
examples we discuss, notably the characterization of the resonance.
We will thus use the ARS algorithm for two purposes: first, when the
resonance equation has no positive integral resonance or when all
resonances are compatible, then the solution is given by a Laurent
expansion; second, when Laurent expansion fails, we use the ARS
algorithm to guess the possible form of the psi-series expansion we
are looking for, and then the proof will be conducted along the same
way we do for $p_n$. Of course, there are also cases for which the
ARS algorithm can be easily justified and the singularity is not
essential (say, by the absolute convergence of the psi-series).

\paragraph{Absolute convergence of the psi-series.} We now prove that
$U(Z)$ converges absolutely in a cut-disk $\mathscr{C}_R$ for some
positive $R>0$.

\begin{prop} \label{prop-ac} For each fixed $c_6$, the psi-series
expansion (\ref{UZ-exp}) converges absolutely for $z$ in the
cut-disk $\mathscr{C}_{(1-\ve)\rho}$ (defined in (\ref{CR})), where
$\ve>0$ is a small number.
\end{prop}
The range $|z-\rho|\le (1-\ve)\rho$ is the best that our approach can
achieve although it seems to hold true, by numerical evidence, up to
$|z-\rho|\le\rho$; in particular, this suggests that the psi-series
expansion be convergent even for $Z=1$ or $z=0$ for $P(z)$.

From this proposition, we see that the solution $P(z)$ can be
analytically continued to at least the region
\[
    \left\{\{ z\,:\, |z|\le \rho+\ve\} \cup
    \{ z\,:\, |z-\rho|\le (1-\ve)\rho\}\right\}
    \setminus [\rho,(2-\ve)\rho]\qquad(\ve>0),
\]
from which we deduce (\ref{pnn}).

To prove Proposition~\ref{prop-ac}, we adopt an approach due to
Hille \cite{Hille76} with some new ingredients; see also
\cite{Hille73}. The resulting proof can then be extended to cover
all the types of DEs we discuss in this paper, whatever their
orders.

\paragraph{Proof of the absolute convergence of the psi-series.
I. Recurrence of $u_k$.} We first rewrite the DE (\ref{bst-de}) for
$P$ into that for $U$, which becomes
\begin{align*}
    \left((1- Z)U'(Z)\right)'=\rho U(Z)^2.
\end{align*}
For convenience, let $U_0=\rho U$. Then
\begin{align*}
    \left((1- Z)U_0'(Z)\right)'=U_0(Z)^2.
\end{align*}
As in \cite{Hille73}, we then convert this DE into a first-order
differential system by introducing an additional function $V_0
:=(1-Z)U_0'(Z)$ as follows.
\begin{align}\label{Bst-de-sys}
    \left\{\begin{array}{l}
        \ds U_0'(Z)=\frac{V_0(Z)}{1-Z},\\
        V_0'(Z)= U_0(Z)^2.
    \end{array}\right.
\end{align}
Let $\tau=\log Z$, $U_0(Z)=\sum_{k\ge0} u_k(\tau) Z^{k-2}$ and
$V_0(Z)=\sum_{k\ge0} v_k(\tau) Z^{k-3}$, where $u_k$ and $v_k$ are
polynomials in $\tau$ of degree at most $\lfloor k/6\rfloor$. Note
that $(\text{d}\tau)/(\text{d}Z)= Z^{-1}$ and $c_0=6/\rho$. From
(\ref{Bst-de-sys}), we derive an infinite system of equations in $k$
($\dot{u}_k := u_k'(\tau)$)
\[
\begin{split}
    \left\{\begin{array}{l}
        \ds \dot{u}_k+(k-2)u_k=v_k+\sum_{0\le j<k}v_j, \\
        \ds \dot{v}_k+(k-3)v_k=12u_k+\sum_{1\le j<k}
        u_ju_{k-j},
    \end{array}\right.
\end{split}\qquad (k\ge 7).
\]
We can further express the above system in terms of matrices as
follows. Let
\[
    \bm{\phi}_k :=\left(\begin{array}{c}
        u_k\\ v_k
    \end{array}\right),\;
    {\mathbf A_k} :=\left(\begin{array}{cc}
        k-2 & -1 \\ -12 &k-3
    \end{array}\right),\;\text{and}\;
    \mathbf{g}_k:=\left(\begin{array}{c}
        \ds\sum_{0\le j<k}v_j \\
        \ds\sum_{1\le j<k}u_ju_{k-j}
    \end{array}\right).
\]
Then, for $k\ge 7$,
\begin{align}\label{Pfsys}
    \dot{\bm{\phi}_k}+{\mathbf{A}_k\bm{\phi}_k} =\mathbf{g}_k,
\end{align}
which can be explicitly solved.
\begin{lmm} For $k\ge 7$, $\bm{\phi_k}$ admits a unique solution
satisfying \label{lmm-phi-k}
\[
    \lim_{\tau\to -\infty}
    \|e^{\mathbf{A}_k\tau}\bm{\phi}_k(\tau)\|=0
\]
of the form
\begin{align}\label{phi-k}
    \bm{\phi}_k(\tau)&=\int_0^\infty
    e^{-x\mathbf{A}_k}\mathbf{g}_k(\tau-x)\dd{x}\nonumber \\
    &=\int^{\infty}_0 \mathbf{P}\,e^{-x\mathbf{D}}
    \mathbf{P}^{-1}\mathbf{g}_k(\tau-x)\dd{x},
\end{align}
where $\mathbf{D}:= \left(\begin{array}{cc} k+1 & 0 \\ 0 & k-6
\end{array}\right)$, $\mathbf{P}=\left(\begin{array}{cc} 1& 1 \\ -3 & 4
\end{array}\right)$ and $\mathbf{P}^{-1}=\left(\begin{array}{cc}
\frac 47& -\frac {1}7 \\ \frac 37 & \frac 17
\end{array}\right)$.
\end{lmm}
\pf The fundamental matrix solution associated with the homogeneous
part of (\ref{Pfsys}) is $e^{\tau\mathbf{A}_k}$, so we can solve
(\ref{Pfsys}) by multiplying it by $e^{x\mathbf{A}_k}$ and then by using
the fact that $u_k(\tau)$ and $v_k(\tau)$ are polynomials in $\tau$,
which gives
\begin{align*}
    &e^{x\mathbf{A}_k}\dot{\bm{\phi}_k}(x) +{\mathbf{A}_k
    e^{x\mathbf{A}_k}\bm{\phi}_k}(x)
    = \frac {\dd{}}{\dd{x}}
    \left(e^{x\mathbf{A}_k}\bm{\phi}_k(x)\right)
    =e^{x\mathbf{A}_k}\mathbf{g}_k(x).
\end{align*}
Integrating both sides from $-\infty$ to $\tau$, we get
\[
    e^{x\mathbf{A}_k}\bm{\phi}_k\biggr|_{-\infty}^\tau
    =e^{\tau\mathbf{A}_k}\bm{\phi}_k(\tau)
    =\int_{-\infty}^\tau e^{x\mathbf{A}_k}\mathbf{g}_k (x)\dd{x},
\]
or
\[
    \bm{\phi}_k(\tau)
    =\int_{-\infty}^\tau e^{(x-\tau)\mathbf{A}_k}\mathbf{g}_k (x)\dd{x}.
\]
The lemma then follows by a change of variables. \qed

\paragraph{Proof of the absolute convergence of the psi-series. II.
An estimate for $u_k$.} To estimate the growth order of $u_k$ and
$v_k$, we now introduce the following norm: for any $\mathbf{x} \in
\mathbb{C}^n$ and any matrix $\left(a_{ij}\right)_{n\times n}$,
\[
    \|\mathbf{x} \|=\max_{1\le j\le n}\{ |x_j|\}, \quad
    \|\left(a_{ij}\right)_{n\times n} \|
    =\max_{1\le j\le n}\left\{\sum_{i}|a_{ij}|\right\}.
\]

With this norm, we then have the inequality
\begin{align}  \nonumber
    \max\{|u_{k}(\tau)|,|v_{k}(\tau)|\}
    &\le \|\bm{\phi}_{k} \| \\
    &\le 5 \int^{\infty}_0 e^{-x(k-6)}
    \max\left\{{\sum_{0\le j<k}|v_j|},
    \sum_{1\le j<k}|u_j||u_{k-j}|\right\}\dd{x}.
    \label{fk-bound}
\end{align}
Now write $z=\rho-r e^{\bm{i}\theta}$, so that $\tau=\log(r/\rho)
+\bm{i}\theta=\xi+\bm{i}\theta$, where $r\le e^{-\ve}\rho$ and
\[
   \mathscr{T}:=\left\{\xi+\bm{i}\theta :
   \xi\in (-\infty,-\ve]\;\mbox{and}\;
   |\theta|\le \pi\right\},
\]
with $|1-\tau|\ge 1+\ve$. We prove by induction that
\begin{align}
    \left\{\begin{array}{l}
        |u_k(\tau)| \le \displaystyle
        \frac{K |1-\tau|^{k-6}}{\sqrt{k+1}},\\
        |v_k(\tau)| \le \displaystyle
        \frac{K|1-\tau|^{k-6}}{\sqrt{k+1}},
    \end{array}\right.  \label{fk-est}
\end{align}
for $k\ge0$ and $\tau\in\mathscr{T}$, where the constant $K>0$ is
easily tuned according to the initial conditions.

Then, by induction hypothesis,
\begin{align*}
    \left|\sum_{0\le j<k} v_j(\tau)\right|
    &\le K\sum_{0\le j<k} \frac{|1-\tau|^{j-6}}
    {\sqrt{j+1}} \\
    &\le \frac{K}{|1-\tau|-1}\,|1-\tau|^{k-6}\\
    &\le \frac{K}{\ve} \,|1-\tau|^{k-6},
\end{align*}
and
\begin{align*}
    \left|\sum_{1\le j<k} u_j(\tau)u_{k-j}(\tau)\right|
    &\le K^2 |1-\tau|^{k-12}
    \sum_{1\le j<k} \frac{1}{\sqrt{(j+1)(k-j+1)}} \\
    &\le K^2 |1-\tau|^{k-12}\int_0^k
    \frac1{\sqrt{x(k-x)}}\mbox{d}x\\
    &= \pi K^2 |1-\tau|^{k-12}.
\end{align*}
Now
\begin{align*}
    \max\{|u_k(\tau)|,|v_k(\tau)|\}
    &\le \Arrowvert \bm{\phi}_k(\tau)\Arrowvert\\
    &=\left\Arrowvert
    \int_{-\infty}^\tau \mathbf{P}e^{(x-\tau)\mathbf{D}}
    \mathbf{P}^{-1} \mathbf{g}_k(x) \mbox{d} x\right\Arrowvert \\
    &= \left\Arrowvert \int_0^\infty \mathbf{P}e^{-x\mathbf{D}}
    \mathbf{P}^{-1} \mathbf{g}_k(\tau-x) \mbox{d} x\right
    \Arrowvert \\
    &\le \Arrowvert \mathbf{P}\Arrowvert
    \Arrowvert\mathbf{P}^{-1}\Arrowvert\int_0^\infty
    e^{-x(k-6)}\Arrowvert \mathbf{g}_k(\tau-x)
    \Arrowvert\mbox{d} x \\
    &\le 5\int_0^\infty e^{-x(k-6)}\max\left\{
    \sum_{0\le j<k} \left|v_j(\tau-x)\right|,\right.\\
    &\left. \qquad \qquad
    \sum_{1\le j<k} \left|u_j(\tau-x)u_{k-j}(\tau-x)
    \right|\right\}\mbox{d}x.
\end{align*}
By choosing $\ve\le 1/(\pi K)$, so that $K/\ve \ge \pi K^2$. We have
\begin{align*}
    |u_{k+6}(\tau)|,|v_{k+6}(\tau)|
    &\le \frac{5K}{\ve} \int_0^\infty
    e^{-xk}  |1-\tau+x|^{k} \mbox{d}x \\
    &\le \frac{5K}{k\ve } |1-\tau|^{k}
    \int_0^\infty e^{-x} \left|1+\frac{x}{k(1-\tau)}
    \right|^k \mbox{d}x.
\end{align*}

Since $|1-\tau|\ge 1+\ve$ for $\tau\in\mathscr{T}$, we see that
\begin{align}
    \int_0^\infty e^{-x} \left|1+\frac{x}{k(1-\tau)}
    \right|^k \mbox{d}x &\le \int_0^\infty e^{-x}
    \left(1+\frac{x}{k|1-\tau|}\right)^k \mbox{d}x\nonumber \\
    &\le \int_0^\infty e^{-x(1-1/|1-\tau|)} \mbox{d}x\nonumber\\
    &= \frac{|1-\tau|}{|1-\tau|-1}\nonumber \\
    &\le \frac{1+\ve}{\ve}. \label{ineq}
\end{align}
It follows that
\begin{align*}
    &\frac{5K}{k\ve}|1-\tau|^{k}
    \int_0^\infty e^{-x} \left|1+\frac{x}{k(1-\tau)}
    \right|^k \mbox{d}x \\
    &\qquad \le \frac{5K(1+\ve)}{k\ve^2} |1-\tau|^{k}\\
    &\qquad \le \frac{K |1-\tau|^k}{\sqrt{k+7}},
\end{align*}
for $k\ge k_0\ge -7+(1+\ve)^2/\ve^4$. This proves the required
estimate.

\paragraph{Proof of the absolute convergence of the psi-series:
an estimate for $U(Z)$.}

From (\ref{fk-est}), we obtain
\begin{align*}
    \rho|U(Z)| &=
    \left|\sum_{k\ge0} u_k(\tau) e^{(k-2)\tau}\right|\\
    &\le Ke^{-2\Re(\tau)} \sum_{k\ge0} \frac{(|1-\tau|)^{k-6}
    e^{k\Re(\tau)}}{\sqrt{k+1}} \\
    &= O\left(e^{-2\Re(\tau)}
    \left(1-|1-\tau|e^{\Re(\tau)}\right)^{-1/2}\right)\\
    &= O(1),
\end{align*}
provided that
\[
    |1-\tau|e^{\Re(\tau)} < 1.
\]
But this implies that ($\Re(\tau) = r/\rho$)
\[
    r <\frac{\rho}{|1-\tau|}\le \frac\rho{1+\ve} \le (1-\ve')\rho.
\]
This proves that the series (\ref{UZ-exp}) is absolutely convergent
for $z\in\mathscr{C}_{(1-\ve)\rho}$. \qed

\paragraph{Numerical approximations to $\rho$ and $c_6$.}
\label{page-na}
As mentioned in Introduction, $P$ is connected to $U$ by choosing a
point in $[\ve\rho,\rho-\ve]$; then the values of $(\rho,c_6)$ are
determined by solving numerically the two equations $P(z_0)=U(Z_0)$
and $P'(z_0)=-\rho U'(Z_0)$, where $Z_0 := 1-z_0$.

For numerical purposes, we can compute the approximate values of
$P(z_0)$ or $P'(z_0)$ by their corresponding truncated series
expansions using, say the first $N$ terms; for example,
$P(z_0)\approx \sum_{j<N}p_j z_0^j$. The number of terms used
depends on the degree of numerical precision we require, and the
remainder $\sum_{j\ge N}p_jz^j$ can be well estimated by using the
asymptotic expansion (\ref{pnn}). More precisely, for large $N$,
\begin{align}\label{P-tail}
    \sum_{j\ge N} p_j z_0^j =
    \frac{6(z_0/\rho)^N}{\rho-z_0}\left(N +
    \frac{3\rho+2z_0}{5(\rho-z_0)}
    +O\left(N^{-4}\right)\right).
\end{align}
Since $z_0<\rho$, the right-hand side can be made arbitrarily small
by choosing $N$ sufficiently large so that the error introduced is
under control.

Similarly, $U(Z)\approx U_M(Z) := \rho^{-1}\sum_{k<M} u_k(\log Z)
Z^{k-2}$ for a sufficiently large $M$ whose choice can be determined
by the desired degree of precision and the upper bound
(\ref{fk-est}).
\begin{align}\label{U-tail}
    \sum_{k\ge M} u_k(\tau_0) e^{(k-2)\tau_0}
    =O\left(M^{-1/2} |1-\tau_0|^M e^{M\Re(\tau_0)}\right),
\end{align}
where $\tau_0 = \log(Z_0)$.

Note that if $z_0$ is too close to zero, then the remainder
(\ref{P-tail}) for $P$ decreases much faster than that
(\ref{U-tail}) for $U$, and if $z_0$ is too close to $\rho$, then
the converse is true. So the best choice for $z_0$ will be the one
that both remainders are asymptotically of the same order. For
practical use, since $p_n$ is easier to compute than $u_k$, we take
$M = \beta N$ for some $\beta\in(0,1)$. Then we solve the equation
\begin{align}\label{best-z0}
    \left(\frac{z_0}{\rho}\right)^{1/\beta}
    = \left|1-\log\left(1-\frac{z_0}{\rho}
    \right)\right|\left(1-\frac{z_0}{\rho}\right),
\end{align}
(which obviously has a unique real solution for $z_0/\rho\in
(1/2,1)$) to find the best $z_0$.

On the other hand, to compute $u_k$, we take the first entry of
$\bm{\phi}_k$ in (\ref{phi-k}) and obtain the recurrence
\begin{align*}
    u_k(\tau)&= \frac {1}7\int_0^\infty
    \left(3e^{-(k-6)x}+4e^{-(k+1)x}\right)
    \left((k-3)u_{k-1}(\tau-x)+u_{k-1}'(\tau-x)\right)\dd x\\
    &\qquad+\frac17\int_0^\infty
    \left(e^{-(k-6)x}-e^{-(k+1)x}\right) \sum_{1\le j<k}
    u_j(\tau-x)u_{k-j}(\tau-x)\dd{x}\\
    &= u_{k-1}(\tau)+\frac {1}7\int_0^\infty
    \left(9e^{-(k-6)x}-16e^{-(k+1)x}\right)
    u_{k-1}(\tau-x)\dd x\\
    &\qquad+\frac17\int_0^\infty
    \left(e^{-(k-6)x}-e^{-(k+1)x}\right) \sum_{1\le j<k}
    u_j(\tau-x)u_{k-j}(\tau-x)\dd{x},
\end{align*}
for $k\ge 7$. All these polynomials $u_k$'s are solvable recursively
starting from the initial values
\[
    u_0=6,\,u_1=-\tfrac{12}{5},\,u_2=-\tfrac7{25},\,
    u_3=-\tfrac{14}{125},\, u_4=-\tfrac{63}{1250},\,
    u_5=-\tfrac{161}{9375},\,
    u_6=c_6-\tfrac{14\tau}{3125},
\]
with the two free parameters $\rho$ and $c_6$. More explicitly, let
$u_k(\tau) := \sum_{0\le s\le \tr{k/6}} u_{k,s} \tau^s$. Then
\begin{align*}
    u_{k,s} &= u_{k-1,s} + \frac1{s!}\sum_{s\le \ell \le
    \tr{(k-1)/6}} u_{k-1,\ell} (-1)^{\ell-s}
    \left(\frac{9\ell!}{7(k-6)^{\ell-s+1}}-\frac{16\ell!}
    {7(k+1)^{\ell-s+1}}\right)\\
    &\; + \frac1{s!}\sum_{\substack{1\le j<k\\
    0\le \ell_1\le\tr{j/6}\\ 0\le \ell_2\le \tr{(k-j)/6}\\
    \ell_1+\ell_2\ge s}} \!\!\!u_{j,\ell_1} u_{k-j,\ell_2}
    (-1)^{\ell_1+\ell_2-s} \left(\frac{(\ell_1+\ell_2-s)!}
    {7(k-6)^{\ell_1+\ell_2-s+1}}-\frac{(\ell_1+\ell_2-s)!}
    {7(k+1)^{\ell_2+\ell_2-s+1}}\right),
\end{align*}
for $0\le s \le \tr{k/6}$.

We finally solve numerically the pair $(\rho,c_6)$ from the two
equations with $\rho\in(3,4)$
\begin{align}
    P_N(z_0)=U_M(Z_0) \quad\text{and} \quad
    P_N'(z_0)=-\rho U_M'(Z_0). \label{z0z1}
\end{align}

Numerical evidence suggests that the series definition for $U(Z)$
and $U'(Z)$ are both convergent for $Z=1$, which means that one
might even use the two equations
\begin{align*}
  U(1)=1 , \quad U'(1)=-\rho,
\end{align*}
to solve for the pair $(\rho,c_6)$. But the convergence is much
slower than taking $z_0$ according to (\ref{best-z0}).

\paragraph{A quantity arising in phylogenetic trees.} Very similar to
the original motivations of studying $p_n$, the following recurrence
\begin{align}\label{qn}
   q_n =\frac 2{(n-1)^2}\sum_{1\le j<n} q_j q_{n-j}
   \qquad (n\ge2),
\end{align}
with $q_1=1$ was introduced in Bryant et al.\ \cite{BMS03} in the course
of analyzing the size of a maximum agreement subtree in two randomly
chosen trees according to the Yule-Harding model. The quantity
serves as an effective bound for the probability that the size of a
common maximum agreement subtree exceeds a certain given value.

Let $p_n := 2q_{n+1}$. Then the recurrence (\ref{qn}) becomes
\[
    p_n = n^{-2} \sum_{0\le j<n} p_j p_{n-1-j}\qquad(n\ge1),
\]
of exactly the same form as (\ref{pn}) but with $p_0=2$. This means
that the DE satisfied by the generating function $P(z) = \sum_n p_n
z^n$ remains the same as (\ref{bst-de}) but the initial condition
differs.

The same psi-series method we used above applies and we obtain 
the asymptotic expansion
\[
    q_n =\rho^{-n} \left(3n -\frac 65 +\frac {168}{3125}n^{-5}
    +\frac {336}{3125}n^{-6}+O(n^{-7}) \right).
\]
with $\rho=1.57042\,87836\,01468\,47580\,40837\dots$.

\section{Probability of equality of random trees}
\label{sec:rts}

The consideration of the equality of two random BSTs can be easily
extended either to more random BSTs or to other variants of BSTs.

\subsection{Equality of $d$ random BSTs}
\label{sec:d-bst}

We extend in this subsection the same psi-series analysis to $d$
random BSTs, $d\ge2$. Surprisingly, the resulting forms of the
asymptotic expansions depends on the parity of $d$.

\paragraph{Recurrence.} The random BST model is as introduced above.
Let $p_n=p_n(d)$ denote the probability that $d$ random BSTs, each
independent of the others, are identical. More precisely, the
probability that $d$ random permutations whose corresponding BSTs
are all the same. Then $p_n$ satisfies the recurrence
\begin{align}\label{pnd-rr}
    p_n = n^{-d} \sum_{0\le j<n} p_j
    p_{n-1-j}\qquad(n\ge1),
\end{align}
with $p_0=1$. Let $P(z) :=\sum_{n\ge 0} p_nz^n$ be the generating
function of $p_n$. Then $P(z)$ satisfies the nonlinear DE of order
$d$
\begin{align} \label{d-bst-equal}
   \left( z\frac {\text{d}}{\text{d}z}\right)^d P(z)=zP(z)^2
\end{align}
with $p_0=1$ and the first $d-1$ values $p_n$ for $1\le n<d$ given
by the recurrence (\ref{pnd-rr}).

\paragraph{The ARS Algorithm.} As in the case of two random BSTs
above, we begin with applying the ARS Algorithm and check first if
there are pseudo-poles and incompatibility.

\begin{itemize}

\item[\ding{182}] Leading order analysis: This part is always easy
for the problems we study in this paper and we obtain, by assuming
$P(z)\sim c_0(1-z/\rho)^{-\alpha}$ and by matching coefficients,
$\alpha=d$ and $c_0=\rho(2d)!/(2 d!)$.

\item[\ding{183}] Resonance analysis: On the other hand, by
collecting the coefficient for the term $c_r(1-z/\rho)^{r-2d}$ in
the resulting expansion for (\ref{d-bst-equal}), we obtain the
polynomial characterizing all possible resonances
\begin{align}\label{d-bst-reson}
    \Phi_d(r)&= \frac{(2d-1-r)!}{(d-1-r)!}
    -\frac {(2d)!}{d!}\\ \nonumber
    &=\left\{\begin{array}{ll}
        (r+1)\phi_d(r),\quad &\mbox{$d$ is odd;}\\
        (r+1)(r-3d)\phi_d(r),\quad &\mbox{$d$ is even,}
    \end{array} \right. \quad d\in\mathbb{N},
\end{align}
where $\phi_d$ is a polynomial of even order and has no real zeroes.
We see that if $d$ is odd, then there is no additional
integer-valued resonance except $-1$ for this case. Thus, the
movable singularity $\rho$ is a pole of order $d$. On the other
hand, if $d$ is even, then there exists an additional, unique,
positive, integer-valued resonance $3d$ for each $d$.

\item[\ding{184}] Incompatibility: We need only consider the case
when $d$ is even. The incompatibility of the resonance at $r=3d$ is
easily checked for each specific $d=2,3,\dots$, but a proof that
$r=3d$ leads to incompatibility for all $d$ is not obvious.
\end{itemize}

\paragraph{The case when $d$ is odd.} From the above quick check
by ARS algorithm, we see that the solution for the DE
(\ref{d-bst-equal}) admits the Laurent series expansion
\[
    \rho P(z) = \frac {(2d)!}{2\cdot d!} \left( Z^{-d}-
    \frac {(3d-2)(d-1)}{2(3d-1)}
     Z^{-d+1}+\sum_{2\le j\le d}c_j Z^{j-d}\right)+\Xi(z),
\]
where $\Xi(z)=\Xi_d(z)$ is analytic at $\rho$.

\paragraph{The case when $d$ is even.} By the above procedure of
ARS algorithm, we anticipate a psi-series expansion for $P(z)$ of
the form
\begin{align}\label{Pdz}
    \rho P(z) = \sum_{j\ge 0} Z^{j-d} \sum_{0\le \ell\le
    \lfloor j/3d\rfloor} c_{j,\ell} (\log Z)^\ell,
\end{align}
where the $c_{j,\ell}$'s are chosen so that the psi-series satisfies
the DE (\ref{d-bst-equal}). In particular, the first few terms read
\[
    \rho P(z)= \frac {(2d)!}{2\cdot d!} Z^{-d}
    -\frac {(3d-2)(d-1)(2d)!}{4(3d-1)d!}
     Z^{1-d}+\sum_{2\le j\le 3d}c_{j,0}Z^{j-d}
     + C_{3d,1} Z^{2d}\log Z+\cdots.
\]
The justification of the psi-series on the right-hand side of
(\ref{Pdz}) follows the same pattern as that for two random BSTs;
see Appendix A1 for details.

In summary, we conclude the following asymptotic estimates, the
drastic change of the error term according to the parity of $d$
unveiling an additional surprise.

\begin{thm}\label{thm-2} The probability that $d\ge2$ randomly chosen
BSTs are all equal satisfies
\begin{align*}
    p_n  &=\rho^{-n-1}\frac {(2d-1)!}{(d-1)!^2}
    \left(n^{d-1}+\frac{(d-1)(2d-1)}{3d-1}\,n^{d-2}
    + \sum_{0\le j\le d-3} C_j n^j\right)\\
    &\qquad+ \left\{\begin{array}{ll}
        O(\rho^{-n}(1-\ve)^n),& \mbox{if $d$ is odd;}\\
        \displaystyle Kn^{-2d-1}\rho^{-n-1}+
        O\left(\rho^{-n}n^{-2d-2}\right), &\mbox{if $d$ is even,}
    \end{array} \right.
\end{align*}
where $\ve>0$, the $C_j$'s are constants, $\rho=\rho_d$
depends on $d$ and $K$ is a constant depending only on $d$.
\end{thm}
More precise asymptotic expansions can be derived, but we content
ourselves with the current form for simplicity of presentation. Is
there any intuitive reason why the asymptotic expansion of
$p_n=p_n(d)$ differs according to the parity of $d$?

\subsection{Equality of two random $m$-ary search trees}
\label{sec:2mary}

The $m$-ary search trees are one of the natural extensions of BSTs
to branching factors $m\ge2$ beyond binary; see \cite{Mahmoud92} for
thorough discussions. Briefly, the first $m-1$ keys are stored in
the root and sorted in increasing order, each of the remaining
$n-m+1$ keys are then directed to one of the $m$ subtrees,
corresponding to the $m$ intervals specified by the $m-1$ sorted
keys, and are constructed recursively by the same procedure.

In the same vein, the probability $q_n$ that two random $m$-ary
search trees are identical is characterized by the following
recurrence ($m\ge2$)
\[
    q_n=\binom{n}{m-1}^{-2}
    \sum_{\substack{j_1+\cdots+j_m=n-m+1\\ j_1,\ldots,j_m\ge0}}
    q_{j_1} \cdots q_{j_m} \qquad(n\ge m-1),
\]
with the initial conditions $q_j=1$, $0\le j\le m-2$. The associated
generating function $Q(z)$ then satisfies the following nonlinear DE
\begin{align} \label{Pz-mary}
    \left(z^{m-1} Q^{(m-1)}(z) \right)^{(m-1)}
    = (m-1)!^2 Q^m(z),
\end{align}
with the initial conditions $Q(z) = 1 + z + \cdots + z^{m-2} +
q_{m-1}z^{m-1}+\cdots$ where $q_j$, $m-1\le j\le 2m-3$, are
determined by the above recurrence.

\begin{itemize}

\item[\ding{182}] Leading order analysis: The simple form $Q(z) \sim
c_0 (1-z/\rho)^{-\alpha}$ leads to $\alpha=-2$ and $\rho c_0=\left(
{(2m-1)!}/{(m-1)!^2}\right)^{1/(m-1)}$.

\item[\ding{183}] Resonance analysis: Again, assuming that $Q(z)
\sim c_0(1-z/\rho)^{-2} + c_r(1-z/\rho)^{-2+r}$, we obtain the
following algebraic equation characterizing all possible resonances
\[
    \prod_{2\le j<2m}(r-j)-\frac {(2m)!}{2}
    =(r+1)(r-(2m+2))\phi_m(r)=0,
\]
where $\phi_m(r)$ is a polynomial of degree $2(m-2)$ and admits
complex-conjugate zeros only. Thus we need to check if the DE
(\ref{Pz-mary}) is compatible at the resonance $r=2m+2$.

\item[\ding{184}] Incompatibility: Similar to the case of $d$ random
BSTs, the resonance $r=2m+2$ is easily checked to be incompatible
for each finite values of $m=2,3,\dots$, but it is far from being
obvious to prove directly the incompatibility for all $m\ge2$.
\end{itemize}

Let $\lambda_m:=\ds \left((2m-1)!/(m-1)!^2\right)^{1/(m-1)}$.
Instead of proving the incompatibility of $r=2m+2$ for all $m\ge2$
and that $\rho$ is not an essential singularity, we prove that the
DE (\ref{Pz-mary}) has the psi-series solution
\[
    U(Z) = \sum_{j\ge0} Z^{j-2} \sum_{0\le \ell \le \lfloor j/(2m+2)
    \rfloor}c_{j,\ell} \log^\ell Z,
\]
which converges absolutely in some cut-region $\mathscr{C}_R$
(defined in (\ref{CR})); see Appendix A1 for details. Then we connect
$Q(z)$ and $U(Z)$ by the same arguments as those used above for two
random BSTs. In this way, we obtain
\begin{align*}
    \rho Q(z) &=\lambda_m Z^{-2} -\frac
    {m\lambda_m}{2m+1} Z^{-1}+\sum_{2\le j\le 2m+2}
    c_{j,0} Z^{j-2} \\
    &\qquad +c_{2m+2,1}Z^{2m}\log Z+
    O\left(  Z^{2m+1}\log  Z\right).
\end{align*}
From this expansion, we then derive the following approximation to
$q_n$.
\begin{thm} The probability $q_n=q_n(m)$ that two random $m$-ary
search trees are equal satisfies the asymptotic approximation
\[
    q_n  = \lambda_m\rho^{-n-1} \left(n +\frac{m+1}{2m+1}
   \right) + K \rho^{-n-1} n^{-2m-1}+O\left(\rho^{-n}n^{-2m-2}\right),
\]
where $\rho=\rho_m$ and $K$ both depend on $m$.
\end{thm}

\begin{table}
\begin{center}
{\setlength{\fboxrule}{2mm}
\begin{tabular}{|c|l|r|}\hline
    $m$ &\hfil $p_n \sim $\hfil & \hfil
    $\lambda_m$\hfil\mbox{} \\ \hline \hline
    $2$ & $\lambda_2\rho_2^{-n-1} \left (n +\frac{3}{5} +
    \frac{56}{3125} \,n^{-5}\right)$ & 6 \\ \hline
    $3$ & $\lambda_3\rho_3^{-n-1}\left( n+\tfrac 47
    +{\frac {6927696}{78236585}}\,n^{-7} \right)$
    & $\sqrt{30}$ \\ \hline
    $4$ & $\lambda_4\rho_4^{-n-1} \left(n+\frac 5{9}
    +\frac {10419284224}{15568564095}\,n^{-9}\right)$
    & $\sqrt[3]{140}$ \\ \hline
    $5$ & $ \lambda_5\rho_5^{-n-1} \left(n+\frac 6{11}
    +\frac{1526061507281984000}{194179984589469879}\,
    n^{-11}\right)$ & $\sqrt[4]{630}$\\ \hline
    $6$ & $\lambda_6 \rho_6^{-n-1}\left(n+\frac 7{13}
    +\frac{132275788517112977050000}
    {942913507718961369877}\,n^{-13}\right)$
    & $\sqrt[5]{2772}$\\ \hline
\end{tabular}}
\end{center}
\caption{The asymptotic approximation to the probability that two
random $m$-ary search trees are equal for $m=2,\dots,6$. All
$O$-terms are omitted.}
\end{table}

As for BSTs, the consideration can be extended to choose $d\ge2$
random $m$-ary search trees, and the resonance equation is given by
\[
    \prod_{0\le j<d(m-1)} (d-r+j) - \frac{m(dm-1)!}{(d-1)!}
    = \frac{\Gamma(d-r+d(m-1))}{\Gamma(d-r)}-
    \frac{m(dm-1)!}{(d-1)!}.
\]
We then deduce that this equation has no positive integral resonance
when $m$ is even and $d$ is odd, and has the positive resonance
$d(m+1)$ for all other cases with $d, m\ge2$. Our approach can be
applied and we obtain an asymptotic approximation to the probability
that $d$ random $m$-ary search trees are equal, \emph{the error
terms beyond the constant term being either exponentially small when
$m$ is even and $d$ is odd or of order $\asymp n^{-dm-1}$ for all
the remaining meaningful cases}.

\subsection{Equality of two random fringe-balanced BSTs}
\label{sec:2fbbst}

Median-of-$(2t+1)$ (or fringe-balanced) BSTs represent yet another
class of extensions of BSTs. The idea is, instead of placing the
first element in the given sequence at the root, which may result in
a less balanced binary tree, we take a small sample of size $2t+1$
and use the median of this sample as the root element, which then
partitions the remaining elements as in the construction of BSTs,
where $t\ge0$. This simple balancing scheme has turned out to be
useful for small $t$, notably for the corresponding quicksort
algorithm. Note the the original BST corresponds to $t=0$.

For the probability model, assume, as in random BSTs, that we are
given a random permutation; then we construct the corresponding
median-of-$(2t+1)$ BST, which is called a \emph{random
median-of-$(2t+1)$ BST}.

Let now $f_n=f_n(t)$ denote the probability that two randomly chosen
permutations lead to an identical median-of-$(2t+1)$ BST. Then $f_n$
satisfies the recurrence
\begin{align} \label{rr-pn-fb}
    f_n = \sum_{t\le j\le n-1-t} \frac{\binom{j}{t}^2
    \binom{n-1-j}{t}^2}{\binom{n}{2t+1}^2}f_j f_{n-1-j}
    \qquad(n\ge 2t+1),
\end{align}
with the initial conditions $f_n=1$ for $0\le n\le 2t$.

Let $F(z):=\sum_{n\ge0} f_n z^n$ denote the generating function of
$f_n$. Then $F(z)$ satisfies the DE
\begin{align} \label{de-m2tp1}
    \left( z^{2t+1} F^{(2t+1)}(z) \right)^{(2t+1)}
    = \frac{(2t+1)!^2}{t!^4} \left(
    \left( z^t F^{(t)}\right)^{(t)}(z) \right)^2,
\end{align}
with the initial conditions $F^{(j)}(0) = j!$, $0\le j \le 2t$, and
$f_j$, $2t+1\le j\le 4t+1$, given by the recurrence
(\ref{rr-pn-fb}).

\begin{itemize}

\item[\ding{182}] Leading order analysis: With the simple form $F(z)
\sim c_0 (1-z/\rho)^{-\alpha}$, we obtain $\alpha=2$ and
\[
   \rho c_0 =\frac {(4t+3)!\,t!^4}{(2t+1)!^4},
\]
for each $t\ge 0$.

\item[\ding{183}] Resonance analysis: Again, assuming that $F(z)
\sim c_0(1-z/\rho)^{-2} + c_r(1-z/\rho)^{-2+r}$, we obtain the
resonance equation
\begin{align*}
    \Phi_t(r)&=\left(\prod_{2\le j\le 2t+1} (r-j)\right)
    \left(\prod_{2t+2\le j \le 4t+3}(r-j) - 2
    \prod_{2t+2\le j \le 4t+3}j\right),
\end{align*}
which can be factored into the form
\[
    (r+1)(r-6t-6)\,\phi_t(r)
    \prod_{2\le j\le 2t+1} (r-j),
\]
where $\phi_t(r)$ has only complex conjugate zeros since the factor
\begin{align*}
    \lefteqn{(r-2t-2)\cdots(r-4t-3)-2(2t+2)\cdots(4t  +3)}\\
    &=(r-2t-2)\cdots(r-4t-3)- (2t+3)\cdots(4t+4)
\end{align*}
never vanishes for $r\in\mathbb{R}\setminus\{-1,6t+6\}$. Thus we get
yet another new pattern for the least positive integer-valued
resonance
\begin{align}
    r= \left\{ \begin{array}{ll}
        6, \quad &t=0, \\
        2, & t\ge 1.
    \end{array} \right .
\end{align}

\item[\ding{184}] Incompatibility: As $t=0$ has already been
addressed in Section~\ref{sec:bst}, we focus on $t\ge1$, which has
the constant resonance $r=2$. A direct check of the incompatibility
is possible for $r=2$ and $t\ge1$; see Appendix A2.

\end{itemize}

The same psi-series method applies and we obtain for $t\ge1$
\begin{align*}
    \rho F(z) &=\frac{(4t+3)!t!^4}{(2t+1)!^4}
    \left( Z^{-2}-\frac{2(t+1)^2}{6t+5} Z^{-1}
    +{\frac { \left( 22\,{t}^{2}+35\,t+14\right)
    \left( t+1 \right) ^{2}t }{ \left( 7\,t+6 \right)
    \left(6\,t+5 \right) ^{2}}} \log  Z \right)\\
    & \mbox{}+O\left( |Z||\log Z|)\right).
\end{align*}

\begin{thm} The probability $f_n$ that two random median-of-$(2t+1)$
BSTs are equal satisfies the asymptotic approximation
\begin{align*}
    f_n
    &=\frac {(4t+3)!t!^4}{(2t+1)!^4}\rho^{-n-1}
    \left( n+\frac{3+2t-2t^2}{6t+5}-
    {\frac {\left( 22\,{t}^{2}+35\,t+14 \right)
    \left( t+1 \right) ^{2}t}{ \left( 7\,t+6 \right)
    \left( 6\,t+5 \right) ^{2}}}n^{-1}\right)\\
    & \qquad\mbox{}+O\left(\rho^{-n}n^{-2}\right),
\end{align*}
for $t\ge1$, where $\rho=\rho_t$ is an effectively computable
constant.
\end{thm}
Note that the expansion also holds when $t=0$ but the $O$-term
becomes $O(n^{-5})$; see (\ref{pnn}). Also more terms can be
computed by the same procedure.

\section{Moments of high orders}
\label{sec:mho}

In addition to the equality of random trees, another rich source
where nonlinear recurrences and differential equations of the same
type as we analyzed above arise is the asymptotics of moments of
high orders.

\subsection{Partial match queries in random quadtrees}

We consider first in this section the cost of partial match queries
in random two-dimensional quadtrees. The expected cost was first
analyzed in \cite{FGPR93} (see also \cite{CH03}) and the limit law
derived in \cite{NR99} under an idealized model where randomness is
preserved throughout the tree.

Let $v=(\sqrt{17}-3)/2$. Then the cost of a random partial match
query in a random two-dimensional quadtree of $n$ nodes tends
(under an idealized model where randomness is preserved for all
subtrees), after normalized by $n^v$, to a limit law $X$ whose
moments satisfy (see \cite{NR99})
\[
   \mathbb{E}(X^m) =  \frac{a_m}{\Gamma(mv+1)} ,
\]
where $a_1:=\Gamma(2v+2)/(2\Gamma(v+1)^2)$ and
\begin{align*}
    a_m = \frac{2}{v(m-1)((m+1)v+3)}
    \sum_{1\le j<m} \binom mj a_j a_{m-j}\qquad
    (m\ge2).
\end{align*}

Then the generating function $A(z) := 1+\sum_{m\ge1}a_mz^m/m!$
satisfies the differential equation
\begin{align}\label{Az}
    v^2 z^2A''(z) +2zA'(z)+2A(z) = 2A^2(z),
\end{align}
with the initial conditions $A(0)=1$ and $A'(0)=a_1$.

The psi-series method we use above can be readily applied with the
resonance $r=6$ and we obtain
\begin{align}\label{Az-ae}
\begin{split}
     A(z) &=3v^2 Z^{-2} + \frac65(9v -5)
     Z^{-1}+\sum_{2\le j\le 7} c_j  Z^{j-2}
     +\frac{117(39v+139)}{43750}\,Z^4\log Z \\
     &\qquad +\frac{468(153v+545)}{109375}\,Z^5\log Z
    +\,O\left( |Z|^6|\log Z|\right),
\end{split}
\end{align}
where the $c_j$'s are unimportant constants. By singularity
analysis (\cite{FO90}), we conclude the following asymptotic
approximation to $a_n/n!$.
\begin{thm} The $m$-th moment of $X$ satisfies for large $m$
\begin{equation}\label{EXm}
\begin{split}
    \mathbb{E}(X^m) &=
    \frac{m!\rho^{-m}}{\Gamma(mv+1)}
    \left(3v^2 m +\frac{9}{5}v -
    \frac{1404(39v+139)}{2185\,m^5} \right. \\
    &\hspace*{4cm} \left.+
    \frac{8424(139v+495)}{21875\,m^6}
    + O\left(m^{-7}\right)\right),
\end{split}
\end{equation}
where $\rho \approx
1.37649\,44410\,57156\,25755\dots$.
\end{thm}
We omit all details as they are very similar to the case of the
equality of two random BSTs.

An interesting implication of our psi-series analysis is that we can
derive an asymptotic expansion for the moment generating function of
$X$
\begin{align}\label{MGF-X}
    \mathbb{E}(e^{Xz})  =
    e^{(z/\rho)^{1/v}}\left(3\left(\frac z\rho\right)^{1/v}
    +\frac95 -\frac{22464}{21875}\left(\frac z\rho\right)^{-5/v}
    +O(|z|^{-6/v})\right),
\end{align}
as $|z|\to\infty$ in the sector $|\arg(z)|\le (v-\ve)\pi/2$. This
is proved by the integral representation
\[
    \mathbb{E}(e^{Xz}) = \frac{1}{2\pi i}
    \int_{\mathscr{H}} e^s s^{-1} A(z/s^v) \dd s,
\]
for a suitable Hankel-type contour, and standard analysis; see
Appendix A3. Such an expansion for the moment generating function is
unusual in the probability literature and implies in turn that
\begin{align}\label{log-tails}
    -\log \mathbb{P}(X>t) \sim (1-v) v^{v/(1-v)}
    (\rho t)^{1/(1-v)},
\end{align}
for large $t$, by an application of Tauberian argument; see Section
4.12 of Bingham et al.\ \cite{BGT87}.

Note that the transformations $z=\xi^{-v }$ and $A(z)=2\xi Z(\xi)$
brings the DE (\ref{Az}) to the standard form of the so-called
\emph{Emden's equation}
\[
    \frac{{\rm d}^2}{{\rm d}\xi^2}Z(\xi) = \xi^{-1}Z^2(\xi).
\]
But it is not exactly solvable; see \cite[\S\ 2.3]{PZ95} or
\cite[\S\ 12.4]{Hille76}.

\subsection{Partial match queries in random relaxed $k$-$d$ trees}

In a similar setting, the cost of a random partial match query in a
random relaxed $k$-$d$ trees (see \cite{DEM98}) tends, after proper
normalization, to the limit law $Y$ whose moments satisfy (see
\cite{MPP00})
\[
    \mathbb{E}(Y^m) = \frac{b_m}{\Gamma(m\beta+1)},
\]
where $\beta:=(-1+\sqrt{9-8s/k})/2$ ($s$ out of the $k$ coordinates
in the query pattern is specified, the other $k-s$ being
``don't-cares"), and
\begin{align*}
    b_m = \frac{\beta+1}{(m-1)((m+1)\beta+1)}\sum_{1\le j<
    m} \binom{m}{j} (j\beta+1)b_jb_{m-j} \qquad(m\ge2),
\end{align*}
with
\[
    b_1 = \frac{2\Gamma(2\beta+2)}{\beta(\beta+1)^2(2\beta+1)
    \Gamma^3(\beta+1)}.
\]

It follows that the generating function $B(z):=1+ \sum_{m\ge1}
b_mz^m/m!$ satisfies the nonlinear differential equation
\begin{align}
    \beta z^2B''(z) + (\beta+1)^2 zB'(z) + (\beta+1) B(z)
    = (\beta+1) B^2(z) +\beta(\beta+1)zB'(z)B(z), \label{Q:Bz}
\end{align}
with the initial conditions $B(0)=1$ and $B'(0)=b_1$.

The psi-series method applies with a resonance at $r=2$ and we
obtain the expansion
\begin{align*}
    B(z) &= \frac2{\beta+1} Z^{-1} +\frac{\beta-1}{\beta}
    +c_2Z+\frac{2(\beta-1)(\beta+2)}{3\beta^2(\beta+1)}
    \,Z\log Z +c_3Z^2\\ &\qquad+ \frac{(\beta-1)(\beta+2)(\beta+3)}
    {3\beta^3(\beta+1)}Z^2\log Z+c_4Z^3
    + O\left(|Z|^3|\log Z|\right),
\end{align*}
from which we deduce an asymptotic approximation to higher order
moments of $Y$.
\begin{thm} The $m$-th moment of the limit law $Y$ satisfies
\begin{align*}
     \mathbb{E}(Y^m) &= \frac{2m!\rho^{-m}}
     {(\beta+1)\Gamma(m\beta+1)}\left(1+
     \frac{(\beta-1)(\beta+2)}{3\beta^2 m^2}
     -\frac{(\beta-1)(\beta+2)}{\beta^3 m^3}
     +O\left(m^{-4}\log m\right)\right),
\end{align*}
as $m\to\infty$, where $\rho$ depends on $\beta$.
\end{thm}
Consequences of this expansion can be derived as those for $X$.

\subsection{Recursive partition structures.}
\label{sec:rps}

In the context of recursive interval splitting, Gnedin and Yakubovich 
\cite{GY06} derived the following
recurrence relation for the $m$-th moment $h_m$ of
certain limit law $W$ (satisfying a fixed-point equation with
Dirichlet distribution as prefactors)
\begin{align}
    h_m &=\frac{\Gamma(d+\omega)}
    {\Gamma(\omega)^2\Gamma(m\lambda+d+\omega)}
    \sum_{0\le j\le m}\binom{m}{j}\Gamma(j\lambda+\omega)
    \Gamma((m-j)\lambda+\omega)h_j\,h_{m-j},
    \label{hm}
\end{align}
for $m\ge2$ with $h_0=h_1=1$, where $\lambda, \omega>0$ ($\lambda$
is referred to as the \emph{Malthusian exponent}) and $d=2,3,\dots$.

\paragraph{The case when $d=2$.}
Consider first the simplest case when $d=2$. In this case, the
generating function
\begin{align}\label{hz}
    h(z) := \sum_{m\ge0} \frac{h_m \Gamma(m\lambda+\omega)}
    {m!\Gamma(\omega)}\,z^m,
\end{align}
satisfies the DE (using the relation
$(\lambda+\omega)(\lambda+\omega+1)=2\omega(\omega+1)$)
\[
    vz^2h''(z)+zh'(z)+h(z)=h^2(z),
\]
which is exactly of the type of problems we have been examining in
this paper (cf.\ (\ref{Az})), where for simplicity
\[
    v := \frac {\lambda^2}{\omega(\omega+1)}.
\]
For this DE, we can apply the psi-series method and obtain
$(Z=1-z/\rho)$
\begin{align*}
    h(z)&=6v Z^{-2} - \frac65(6v-1) Z^{-1}
    + \sum_{2\le j\le 6} c_j Z^{j-2} +
    K Z^4 \log Z +O\left(|Z|^5|\log Z|\right),
\end{align*}
where
\begin{align*}
    K := \frac{(v-1)^2(v-6)(6v-1)(2v+3)(3v+2)}{43750v^5}.
\end{align*}
Consequently, we deduce the asymptotic expansion for the moments of
$W$
\[
    h_m =\frac{6m!\Gamma(\omega)\rho^{-m}}
    {\Gamma(m\lambda+\omega)}
    \left(vm-\frac {v-1}5-4Km^{-5}
    +O \left(m^{-6} \right)\right),
\]
for large $m$.

\paragraph{The case when $d\ge2$.}
From the recurrence (\ref{hm}), the generating function $h(y)$
(defined as in (\ref{hz})) satisfies the DE
\[
    y^{1-\omega}\frac{\dd{}^d}{\dd{y}^d}
    \left(h(y^\lambda)y^{d+\omega-1}\right)
    =\rising{\omega}{d} h(y^\lambda)^2,
\]
where $\rising{\omega}{d} = \omega\cdots(\omega+d-1)$ denotes the
rising factorial; see \cite{GY06}. The DE is however less
manageable. We rewrite it as follows. Let $z=y^\lambda$ and
$H(z)=z^\kappa h(z)$, where $\kappa:=(d+\omega-1)/\lambda$. Note
that the Malthusian exponent $\lambda$ satisfies the relation
\[
    \frac{\rising{\omega}{d}}
    {\rising{(\lambda+\omega)}{d}} =\frac12.
\]
Then the function $H(z)$ satisfies the DE
\begin{align}\label{gnedin-r}
    \lambda\theta(\lambda\theta-1)\cdots(\lambda\theta-d+1)H(z)
    =z^{-\kappa}\rising{\omega}{d}H(z)^2,
\end{align}
where the differential operator $\theta$ is defined as
$\theta:=z(\text{d}/\text{d}z)$.

The leading order analysis and the resonance analysis give the
dominant exponent $-d$ and the resonance equation is exactly the
same as (\ref{d-bst-reson}) for all $d\ge2$, namely,
$\rising{(d-r)}{d} -\rising{(d+1)}{d}$. It follows that we have the
same asymptotic pattern for $H$ as the case of $d$ random BSTs.

\paragraph{The case when $d$ is odd.}
The movable singularity $\rho$ is a pole of order $d$ and the
solution $H(z)$ admits the Laurent expansion
\begin{align*}
    \rho^{-\kappa} H(z)
    &= \frac {(2d)!\lambda^d}
    {2\cdot d!\rising{\omega}{d}}
    \sum_{0\le j\le d} c_j Z^{j-d} + \Xi_1(z),
\end{align*}
where
\begin{align}\label{c0c1}
    c_0=1,\quad c_1 = -\frac d2-\frac{(4d-2)\omega
    +(d-1)(5d-2)}{2(3d-1)\lambda},
\end{align}
and $\Xi_1(z)$ is an analytic function at $z=\rho$.

\paragraph{The case when $d$ is even.}
In this case, since the resonance equation (\ref{d-bst-reson})
possesses the unique positive integral resonance $3d$, we see that
$z=\rho$ is a pseudo-pole and the psi-series solution to
(\ref{gnedin-r}) has the form
\begin{align*}
    \rho^{-\kappa}H(z)
    &= \sum_{j\ge 0} Z^{j-d} \sum_{0\le \ell\le \lfloor j/3d
    \rfloor}c_{j,\ell} (\log Z)^\ell \\
    &= \frac {(2d)!\lambda^d}
    {2\cdot d!\rising{\omega}{d}}
    \sum_{0\le j\le 3d} c_j Z^{j-d} +KZ^{2d}\log Z
    +O\left(|Z|^{2d+1}|\log Z|\right),
\end{align*}
where, in particular, $c_0$ and $c_1$ are given as in (\ref{c0c1}),
and $K$ is a constant dependent on $\lambda$ and $\omega$.

\paragraph{Expansions for $h$.}
It is not difficult to verify that $h(z)$ and $H(z)$ have the same
dominant singularity $\rho$, dominant exponent $-d$, and the
dominant resonance $3d$. Now by the relation between $h(z)$ and
$H(Z)$: $h(z)=(1- Z)^{-\kappa}\rho^{-\kappa}\,H(z)$, we obtain
\[
    h(z) = \frac {(2d)!\lambda^d}
    {2\cdot d!\rising{\omega}{d}}
    \times\left\{\begin{array}{ll}
        \ds\sum_{0\le j\le d} c_j' Z^{j-d} +\Xi_2(z),&
        \text{if $d$ is odd};\\
        \ds\begin{array}{l}
        \ds \sum_{0\le j\le 3d} c_j'Z^{j-d} +
        K'Z^{2d}\log Z \\
        \;\;+ O\left(|Z|^{2d+1}|\log Z|\right)
        \end{array},&
        \text{if $d$ is even},
    \end{array}\right.
\]
where $c_0'=1$,
\[
    c_1' = \frac d{2}\left(\frac {d+2\omega-1}
    {(3d-1)\lambda}-1\right),
\]
and $\Xi_2$ is analytic at $z=\rho$.

\paragraph{Asymptotics of the moments.} From the expansions
we derived and a similar analysis as for $d$ random BSTs, we can now
conclude the following asymptotic approximations to the limit law
$W$.

\begin{thm} The $m$-th moment $h_m$ of $W$ satisfies
\begin{align*}
    h_m &= \frac{(2d)!\Gamma(\omega)^2\lambda^dm!\rho^{-m}}
    {2\cdot d!(d-1)!\Gamma(\omega+d)\Gamma(m\lambda+\omega)}
    \sum_{0\le j\le d} C_j m^{d-1-j}\\ &\qquad+
    \left\{\begin{array}{ll}
        O((1-\ve)^m),& \mbox{if $d$ is odd;}\\
        \displaystyle C m^{-2d-1}+
        O\left(m^{-2d-2}\right), &\mbox{if $d$ is even,}
    \end{array} \right.
\end{align*}
where $\ve\in(0,1)$, the $C_j$ are constants with $C_0=1$ and
\[
   C_1 =\binom{d}{2}\frac {d+2\omega-1}{(3d-1)\lambda},
\]
and $\rho, C$ are constants depending on $d, \lambda, \omega$.
\end{thm}

\subsection{An Ansatz solution in Boltzmann equations}
\label{sec:boltzmann}

The following sequence $t_n$ arose in the analysis (see \cite{BC81})
of exact
solutions of the Tjon-Wu representation of Boltzmann equations
(which represent the major cornerstone of kinetic theory in
statistical mechanics). Let $\nu$ be a positive integer. The
sequence $t_n$ is defined recursively as
\begin{align}\label{boltz-rr}
    \left(\frac {\nu(\nu+1)}{\nu+2}n(n-1)-(n+1)
    \right) t_n =-\sum_{0\le j\le n}  t_j  t_{n-j}
    \qquad(n\ge2),
\end{align}
with $ t_0= t_1=1$. This recurrence translates into the following DE
for the generating function $T(z):= \sum_{n\ge0} t_n z^n$
\begin{align}\label{Boltzmann-EqnN}
    \frac {\nu(\nu+1)}{\nu+2} z^2T''(z)- zT'(z)
    -T(z)\left(1-T(z)\right)=0,
\end{align}
with the initial conditions $T (0)=T'(0)=1$.

Straightforward computations as above give $-2$ as the dominant
exponent for the dominant term of $T(z)$ and $(r+1)(r-6)$ as the
resonance equation for each $\nu=1,2,\dots$. Interestingly, for the
resonance $r=6$, the two special cases $\nu=1, 2$ do not lead to
incompatible system of equations, in contrast to all higher values
of $\nu$. This is very different from the cases we have been dealing
with up to now. According to the ARS method, the cases when $\nu=1,
2$ admit the {\em Painlev\'{e} property}\;\cite[\S 1.2,
Definition~1.1]{CM08} and have solutions in terms of Laurent
expansion with two free parameters; in other words, they are
\emph{integrable}, and we will derive closed-form solutions for
them. The remaining cases when $\nu\ge3$ have psi-series solutions.

\paragraph{Exactly solvable (integrable) case : $\nu=1$.}  We
start with the case $\nu=1$. Consider the transformations
$T(z)=1-\zeta V(\zeta)$ and $z=-\zeta$. Note that, by this
transform, the coefficients $[\zeta^n]V(\zeta)$ are positive and the
transformed DE  (after multiplying $V'(\zeta)$) becomes
\[
    \frac 13\zeta^2 \frac {\dd{}}{\dd{\zeta}}
    \left ( \zeta \left(\frac {\dd{V}}
    {\dd{\zeta}}\right)^2-V(\zeta)^3\right)=0;
\]
or equivalently,
\begin{align}\label{boltz-V-sol}
     \sqrt{\zeta}\;\frac {\dd{V}}{\dd{\zeta}}
     =\sqrt{V(\zeta)^3-1}, \quad V(0)=1.
\end{align}
By the relation between $T(z)$ and $V(\zeta)$, we deduce that
$V(0)=1$ and $V'(0)=3$. Then (\ref{boltz-V-sol}) is solved as
\begin{align}\label{imp-eq}
    2\sqrt{\zeta}=\int_{1}^{V(\zeta)} \frac{\dd{x}}{\sqrt{x^3-1}}.
\end{align}

Let
\[
    2\sqrt{\zeta_\infty} = \int_1^{\infty} \frac{\dd{x}}{\sqrt{x^3-1}}
    \approx 2.42865\,06478\,87581\,61181\dots,
\]
or $\zeta_\infty\approx 1.47458\,59923\,71192\,48035\dots$.
Obviously $V(\zeta)\to\infty$ as $\zeta\to\zeta_\infty$. Let $\Delta
:= 2(\sqrt{\zeta_\infty}-\sqrt{\zeta})$. Then (\ref{imp-eq}) can be
written as
\[
    \Delta = \int_{V(\zeta)}^\infty \frac{\dd x}{\sqrt{x^3-1}}.
\]
Since $V(\zeta)\to\infty$ as $\zeta\to\zeta_\infty$, we deduce that
\[
    \Delta = 2 V(\zeta)^{-1/2} +\frac16 V(\zeta)^{-7/2} +\frac3{52}
    V(\zeta)^{-13/2}+\frac{5}{152} V(\zeta)^{-19/2} +
    \text{smaller order terms}.
\]
Consequently, by inverting the series (justified by analyticity and
standard arguments), we obtain
\[
    V(\zeta) = 4\Delta^{-2} +\frac{\Delta^4}{112}+
    \frac{\Delta^{10}}{652288}+\frac{\Delta^{16}}{5552275456}
    +\text{smaller order terms}.
\]

Finally, let $\rho:=-\zeta_\infty$ and we obtain
\begin{align*}
    t_n &=[z^n] T(z) = (-1)^{n-1}[\zeta^{n-1}] V(\zeta)
    =\frac{(-1)^{n-1}}{2\pi i} \oint_{|\zeta|=c<\zeta_\infty}
    \zeta^{-n} V(\zeta) \dd{\zeta} \\
    &= \frac{2(-1)^{n-1}}{2\pi i}
    \oint_{|y|=c'<\sqrt{\zeta_\infty}}
    y^{-2n+1} V(y^2) \dd{y} \\
    &\sim 8(-1)^{n-1}[y^{2n-2}]
    \left(2\sqrt{\zeta_\infty}-2y\right)^{-2}\\
    &= (-1)^{n-1}(4n-2) \zeta_\infty^{-n}\\
    &=2(-1)^{n-1}(2n-1)|\rho|^{-n},
\end{align*}
the errors omitted being exponentially smaller.

\paragraph{Exactly solvable (integrable) case : $\nu=2$.}
The case when $\nu=2$ is similar. We now adopt the transformations
$T(z)=1-\zeta^{2} L(\zeta)$ and $z=-\zeta^3$. Then the DE
(\ref{boltz-rr}) becomes
\[
    \frac{\dd{}^2}{\dd{\zeta^2}}L(\zeta)-6L(\zeta)^2=0 \iff
    \frac{\dd{}}{\dd{\zeta}}\left( \frac
    12\left(\frac{\dd{L}}{\dd{\zeta}}\right)^2-2L(\zeta)^3\right)=0,
\]
with the initial values $L(0)=0$ and $L'(0)=1$. Thus, the solution
is given by
\begin{align}\label{imp-eq2}
   \zeta=\int_0^{L(\zeta)}\frac{\dd{x}}{\sqrt{1+4x^3}}.
\end{align}
Let $\zeta_\infty$ denote the dominant singularity of $L(\zeta)$.
Then
\[
    \zeta_\infty=\int_0^\infty \frac{\dd x}{\sqrt{4x^3+1}}
    =\frac{2^{1/3}}{6}{\rm Beta}\left(\frac16,\frac 13\right)
    \thickapprox 1.76663\,87502\,85449\,95731\dots.
\]
Thus the dominant singularity of $T(z)$ when $\nu=2$ is
\[
    \rho=-\zeta_\infty^3
    = -\frac {1}{108}{\rm Beta} \left(\frac16,\frac13 \right)^3
    \thickapprox -5.51370\,15767\,10567\,75506\dots.
\]
Furthermore, from (\ref{imp-eq2}), we have
\[
    \Delta := \zeta_\infty-\zeta
    = \int_{L(\zeta)}^\infty \frac{\dd x}{\sqrt{4x^3+1}},
\]
and, by the same procedure as above,
\begin{align*}
     \Delta = L(\zeta)^{-1/2}-\frac1{56} L(\zeta)^{-7/2} +
     \frac3{1664}L(\zeta)^{-13/2}
     -\frac{5}{19456} L(\zeta)^{-19/2} +\text{smaller order terms},
\end{align*}
for $\zeta\sim\zeta_\infty^{-}$. By inverting the expansion
\[
    L(\zeta) = \Delta^{-2} -\frac{\Delta^4}{28}
    +\frac{\Delta^{10}}{10192} -
    \frac{\Delta^{16}}{5422144}
    +\frac{3\Delta^{22}}{9868302080}-\text{smaller order terms}.
\]
Accordingly,
\begin{align*}
     t_n &=[z^n]T(z) = \frac1{2\pi i} \oint_{|z|=c<|\rho|}
     z^{-n-1} T(z) \dd{z}\\
     &= \frac{3(-1)^n}{2\pi i} \oint_{|\zeta|=c'<|\rho|^{1/3}}
     \zeta^{-3n-1} T(-\zeta^3) \dd{\zeta} \\
     &= \frac{3(-1)^{n-1}}{2\pi i} \oint_{|\zeta|=c'<\zeta_\infty}
     \zeta^{-3n+1} L(\zeta) \dd{\zeta}
     = 3(-1)^{n-1} [\zeta^{3n-2}] L(\zeta) \\
     &\sim 3(-1)^{n-1} [u^{3n-2}]
     \left(\zeta_\infty-\zeta\right)^{-2} \\
     &= 3(-1)^{n-1} (3n-1)\zeta_\infty^{-3n}\\
     &= 3(-1)^{n-1} (3n-1)|\rho|^{-n}.
\end{align*}

Note that we can use the transforms $z=\zeta^2$ and
$T(z)=1-V(\zeta)\zeta^2$ to convert the DE for $\nu=1$ to a DE of
same type (differing only by a constant) as the case for $\nu=2$.
Also both solutions can be expressed in terms of Weierstrass ${\wp}$
functions.

\paragraph{The rest cases : $\nu\ge 3$.} Unlike the preceding two
cases, the rest $\nu$'s no longer lead  to DEs that are solvable by
{\em quadrature}\footnote{A DE is said to be solvable by quadrature
if its solution can be expressed in terms of one or more
integrations.}. Due to incompatibility, we apply again the
psi-series method. Because of the negative sign on the right-hand
side of (\ref{boltz-rr}), we consider the transform $z=-\zeta$ and
$T(z)=1-\zeta V (\zeta)$. Then
\[
   t_n =\left[z^n\right]T (z)
  =(-1)^{n-1}\left[\zeta^{n-1}\right] V (\zeta),
\]
and (\ref{Boltzmann-EqnN}) is translated into
\[
    \frac {\nu(\nu+1)}{\nu+2}\zeta V''(\zeta)
    +\frac {2\nu^2+\nu-2}{\nu+2}V'(\zeta)
    -V(\zeta)^2=0.
\]
Let now $ Z=1-\zeta/\rho$, where $\rho>0$ is the dominant
singularity of $V$ (having all Taylor coefficients positive).
Then we deduce the psi-series expansion for $V$
\begin{align*}
    \rho V (\zeta) &=\frac {6\nu(\nu+1)}{\nu+2} Z^{-2}
    -\frac {6(\nu^2+2\nu+2)}{5(\nu+2)} Z^{-1}
    +\sum_{0\le j\le 5} c_j  Z^j \\
    &\qquad+ K Z^4\log  Z+O\left( |Z|^5|\log Z|\right),
\end{align*}
where
\[
    K := -\frac {(\nu-1)( \nu-2)( \nu+3 )( \nu+4)
    ( 2\nu+1)(2\nu+3)(3\nu+2)(3\nu+4)
    \left(\nu^2+2\nu+2\right)^2}{43750\nu^5(\nu+1)^5
    (\nu+2)} .
\]
This, together with the approximations we derived for $t_n$ in the
two cases $\nu=1, 2$, implies the following asymptotics of $t_n$.

\begin{thm} The sequence $t_n$ satisfies the asymptotic expansion
\begin{align}\label{Boltzmann-dnN-asymG}
\begin{split}
    (-1)^{n-1} t_n
    = \rho^{-n}\Biggl(
    \frac {6\nu( \nu+1)}{\nu+2} n-\frac{6(\nu^2+2\nu+2)}
    {5(\nu+2)}  +\left\{\begin{array}{ll}
        O((1-\ve)^n),& \text{if }\nu=1, 2;\\
        24K n^{-5} + O(n^{-6}), &\text{if } \nu\ge3.
    \end{array}\right.\Biggr)
\end{split}
\end{align}
\end{thm}
Note that $K=0$ when $\nu=1, 2$.

\section{Conclusions}

Through the examples we studied in this paper, we see that the
psi-series method is a powerful approach to handling nonlinear DEs
and yields several surprising results, notably asymptotic expansions
with the first few terms missing. While psi-series have long been
used in many branches of mathematics and physics, little attention
has been paid to the corresponding asymptotics of the
coefficients. Also the procedure we adapted and improved from Hille
for proving the absolute convergence of psi-series is of certain
generality and can be applied to other problems of similar nature.

Another feature of the recurrences we studied in this paper is that
they are very sensible to small variations, the example of $d$
random BSTs being typical. Note first that the recurrence
(\ref{pnd-rr}) with $d=0$ yields the well-known Catalan numbers and
the case $d=1$ gives rise to the trivial sequence $p_n=1$. The case
$d=1$ in a more general form was studied by Wright \cite{Wright80};
see also Cooper \cite{Cooper47} for a study of $p_n$ for real
$k\ge0$.

We now compare the recurrence (\ref{pnd-rr}) with the following
one by defining $p_1=1$ and
\[
    p_n = n^{-d} \sum_{1\le j\le n-1} p_j p_{n-j}\qquad(n\ge2).
\]
While the case $d=0$ still yields the Catalan numbers with their 
generating function satisfying
\[
    P(z)-z=P^2(z),
\]
the case $d=1$ becomes a nonlinear differential equation of Riccati
type
\[
    zP'(z)-z = P^2(z),\qquad P(0)=0,
\]
which can still be explicitly solved $P(z)=z^{1/2}J_1(2z^{1/2})
/J_0(2z^{1/2})$, where $J_\nu(z)$'s are Bessel functions (see \cite{Ince26}). 
The case
$d=2$ is again of Emden-Fowler type and can be solved asymptotically
by psi-series method as well as the remaining cases $d\ge3$.

See \cite{Cooper47,FGM97,FPV98,Kleitman70,SW79,Wright80}  
and the references therein for some quadratic
recurrences of the above ``Faltung'' type. More examples can be
found in the recent papers \cite{BGR08,BOGRW10}.

\def\theequation{A.\arabic{equation}}
\setcounter{equation}{0}

\section*{Appendix}

\subsection*{A1. Proof of the absolute convergence of psi-series}

In this Appendix, we group the details of the proof of the absolute
convergence of the psi-series arising in the three cases: $d$ random
BSTs, two random $m$-ary search trees, and two random
median-of-$(2t+1)$ BSTs. We first describe briefly the general
pattern of the proof and then provide more details for each case.

Our proof begins with rewriting the original DE in $z$ into a
system of linear DEs in $Z=1-z/\rho$ of the form
\begin{align}\label{A-DE-sys}
    \frac {\dd{}}{\dd{Z}}\mathbf{U}(Z)
    ={\mathcal X}(Z,\mathbf{U}), \quad
    \mathbf{U}(Z)=\left(\begin{array}{c}
        U_{1}(Z) \\
        \vdots \\
        U_s(Z)
    \end{array}\right) ,
\end{align}
where $s\in\{d,2(m-1),4t+2\}$. Here $U_j(Z)=\sum_{k\ge 0}
u_k^{[j]}(\tau) Z^{-\alpha+k-j+1}$, where $\alpha$ is the leading
order, $\tau=\log Z$ and ${\mathcal X}
:\mathbb{C}^{s+1}\mapsto\mathbb{C}^{s}$. Then we derive the infinite
system of linear DEs satisfied by the $u_k^{[j]}$'s
\begin{align}\label{A-infinte-sys}
    \dot{\bm{\phi}_k}+{\mathbf{A}_k\bm{\phi}_k}
    =\mathbf{g}_k,\quad \bm{\phi}_k
    =\left( \begin{array}{c}
        u_k^{[1]} \\
        \vdots \\
        u_k^{[s]}
    \end{array}\right),
\end{align}
where ${\mathbf{A}_k}=k\mathbf{I}_{s\times s}-\mathbf{M}$ and
$\mathbf{M}\in\mathbb{C}^{s\times s}$ are $s\times s$ matrices.

In terms of such an infinite system, an upper bound for all
$u_k^{[j]}$ (in particular, for $u_k^{[1]}$) is of the form
\begin{align*}
    \left|u_k^{[1]}(\tau)\right|
    \le K \psi(k) |1-\tau|^{k-c(s)},
\end{align*}
for $\tau\in\mathscr{T}$
\begin{align}\label{T-region}
    \mathscr{T}:=\left\{\xi+\bm{i}\theta :
    \xi\in(-\infty,-\varepsilon]
    \;\;\mbox{and}\;\;|\theta|\le\pi \right\},
\end{align}
with $|1-\tau|\ge 1+\varepsilon$,  where $K$ is a constant and
$\psi(k), c(s)$ depend on the problem in question. Then the absolute
convergence can be justified.

An additional common and interesting feature this approach brings is
that the resonance equation will be seen to be equal to
$\det(r\mathbf{I}_{s\times s}-\mathbf{M})$. We will explain this in more
details.

The following relations are useful in converting our DEs in $z$ into
those in $Z$ ($\mathbb{D}=\text{d}/\text{d} z$).
\[
    z=\rho(1-Z),\quad z\mathbb{D}=-(1-Z)\frac{\dd{}}{\dd{Z}},
    \quad z^j\mathbb{D}^j=(-1)^j(1-Z)^j\frac{\dd{}^j}{\dd{Z}^j}.
\]

\paragraph{Equality of $d$ random BSTs.}
The corresponding system (\ref{A-DE-sys}) for (\ref{d-bst-equal}) is
\[
  \left\{\begin{array}{ll}
          \ds U_j'(Z)=\frac {U_{j+1}(Z)}{1-Z},
          \quad& 1\le j<d, \\ &\mbox{} \\
          \ds U_d'(Z)=(-1)^d \rho U_1(Z)^2. &
         \end{array}\right.
\]
The associated coefficient matrices $\mathbf{A}_k$ and $\mathbf{g}_k$
in (\ref{A-infinte-sys}), $k\ge 3d+1$,  are given by
\[
     \mathbf{A}_k=k\mathbf{I}_{d\times d}-\mathbf{M}, \quad \mathbf{M}=\left(
     \begin{array}{ccccc}
             d & 1 & 0 & \cdots & 0 \\
             0 & d+1 & 1 & \cdots & 0 \\
             \vdots & \ddots & \ddots & \ddots & \vdots \\
              &  & 0 & 2d-2 & 1 \\
             (-1)^{d-1}\frac{(2d)!}{d!}  & 0 & \cdots & 0 & 2d-1 \\
      \end{array}\right).
\]
and
\[
    \mathbf{g}_k= \left(
    \begin{array}{c}
        \ds\sum_{0\le \ell <k}u_\ell^{[2]}(\tau) \\
        \vdots \\
        \ds \sum_{0\le \ell <k}u_\ell^{[d]}(\tau) \\
        \ds (-1)^d\rho\sum_{1\le \ell <k}
         u_\ell^{[1]}(\tau)u_{k-\ell}^{[1]}(\tau) \\
    \end{array}\right).
\]
Due to the existence of complex-conjugate roots, we can find a
$d\times d$ matrix $\mathbf{P}$ with entries $\mathbf{P}_{ij}\in\mathbb{C}$
such that
\[
    \mathbf{P}\mathbf{A}_k  \mathbf{P}^{-1}
    =\left ( \begin{array}{ccccccc}
        k+1 & 0&\cdots & & &\cdots &0\\
        0&k-3d & & & & &\vdots\\
        & & &k-r_3 & & &\\
        & & & &k-r_4 & &\vdots\\
        \vdots& & & & &\ddots& 0\\
        0&\cdots & & &\cdots & 0&r_d
    \end{array}\right),
\]
for $k\in\mathbb{N}$. By the same norm and same arguments used for
two random BSTs, we derive the inequality
($C_d:=\|\mathbf{P}\|\|\mathbf{P}^{-1}\|$)
\begin{align}\label{appendix-deq-bound}
\begin{split}
    \max_{1\le j\le d}\left\{\left|u_k^{[j]}(\tau)\right|\right\}
    &\le \|\bm{\phi}_k\| \\
    &\le C_d \int^{\infty}_0 e^{-(k-3d)x}\max\left(
    \begin{array}{c}
    \ds\sum_{0\le \ell<k}\left|u_\ell^{[1]}\right|,\ldots,
    \ds\sum_{0\le \ell<k} \left|u_\ell^{[d]}\right|\\
    \ds\rho\sum_{1\le \ell <k}
    \left|u_{\ell}^{[1]}u_{k-\ell}^{[1]}
    \right|\end{array}
    \right)\dd{x}.
\end{split}
\end{align}

Again, by same the arguments used to prove (\ref{fk-est}), we have,
\begin{align*}
    \left|u_k^{[j]}(\tau)\right|
    \le K(1+k)^{-1/2}|1-\tau|^{k-3d}\qquad
    (1\le j\le d, k\ge 0),
\end{align*}
for $\tau\in\mathscr{T}$.

\paragraph{The resonance polynomial equals 
$\det(r\mathbf{I}_{d\times d}-\mathbf{M})$.}
Direct calculations give the determinant
\[
    \det\left(r\mathbf{I}_{d\times d}-\mathbf{M}\right)
    =\frac {(2d-1-r)!}{(d-1-r)!}-\frac{(2d)!}{d!},
\]
which is nothing but the resonance polynomial (\ref{d-bst-reson}).

The reason that the two polynomials are equal is as follows.
The distinction between Laurent
expansion and the psi-series expansion depends crucially either on
the existence of positive integer resonance or on whether a relation
such as (\ref{frob-coeff}) holds for all $k$. This is equivalent to
asking whether the linear system $\mathbf{A}_k\bm{\phi}_k =\mathbf{g}_k$ is
solvable or not for all $k$. If the system (\ref{A-infinte-sys})
$\mathbf{A}_k\bm{\phi}_k=\mathbf{g}_k$ is solvable  under the condition
$\det\mathbf{A}_k\not=0$ for all $k$, then by the uniqueness of the
solution of (\ref{A-infinte-sys}), the solution vectors
$\bm{\phi}_k$'s are constant vectors (independent of $\tau$) and in
turn, the series solution $U_1(Z)=\sum_{k\ge 0} u_k^{[1]}
Z^{-d+k-j+1}$ is eventually a Laurent's series. On the other hand,
if $\det\mathbf{A}_{k_0}\not=0$ fails to hold for some $k_0$, then we
have the following two cases.
\begin{itemize}
\item[---] The linear system $\mathbf{A}_{k_0}\bm{\phi} =\mathbf{g}_{k_0}$ has
a solution depending on the $d-{\rm
rank}\left(\mathbf{A}_{k_0}\right)$ free parameters, and all the rest
constant coefficient vectors $\bm{\phi}_k$ depend on at least these
parameters.

\item[---] The linear system is inconsistent. Hence it can no longer
provide a solution to (\ref{A-infinte-sys}). The real solution
should be solved from (\ref{A-infinte-sys}) instead and then all the
vector functions $\bm{\phi}_k(\tau)$, $k\ge k_0$, depend on
$\tau$. Moreover, the resulting solution $U_1(Z)=Z^{-d}\sum_{k\ge 0}
u_k^{[1]} Z^{k-j+1}$ is indeed a psi-series.
\end{itemize}

In particular, we see that the characteristic polynomial 
$\det(r \mathbf{I}_{d\times d}- \mathbf{M})$ is the same as the
polynomial (\ref{d-bst-reson}) that determines all the possible
resonances.

\paragraph{Equality of two random $m$-ary search trees.}
The transformed first-order differential system in terms of $Z$ for
(\ref{Pz-mary}) now has the form
\begin{eqnarray*}\left\{
\begin{array}{lcl}
    U_1'(Z)=U_{2}(Z) & \quad &1\le j\le m-2, \\
    U_{m-1}'(Z)=(1-Z)^{-(m-1)}U_{m}(Z),&  & \\
    U_{m-1+j}'(Z)=U_{m+j}(Z), & &1\le j\le m-2,  \\
    U_{2m-2}'(Z)=(m-1)!^2\rho^{m-1}U_1(Z)^m. & &
\end{array}\right.
\end{eqnarray*}
So that the corresponding infinite system (\ref{A-infinte-sys}) has
the coefficient matrix $\mathbf{A}_k=k\mathbf{I}_{2(m-1)\times2(m-1)}-\mathbf{M}$,
where
\[
    \mathbf{M}=\left(\begin{array}{crrrrrcccr}
         2 & 1 & 0 & \cdots &&&& \cdots & 0 \\
         0 & 3 & 1 & 0 &&&&& \vdots \\
         \vdots & 0 & 4 & 1 & \ddots &&& \\
         && \ddots & \ddots & \ddots &&&& \\
         &&&& m & 1 &&& \\
         &&&&& m+1 & 1 &&\vdots \\
         \vdots &&&&& \ddots & \ddots & \ddots & 0 \\
         0 &&&&&& 0 & 2m-2 & 1 \\
         m(2m-1)! & 0 & \cdots &&&& \cdots & 0 & 2m-1
     \end{array}\right),
\]
and the vector-valued function
\[
    \mathbf{g}_k=\left(\begin{array}{c}
        0 \\ \vdots \\ 0 \\
        \ds\sum_{0\le j<k} \binom{m-2+k-j}{k-j}u_j^{[m]} \\
        0 \\ \vdots \\ 0 \\
        \ds \rho^{m-1}(m-1)!^2
        \sum_{\substack{i_1+i_2+\cdots+i_m=k\\0\le i_j<k }}
        u_{i_1}^{[1]}u_{i_2}^{[1]}\cdots u_{i_m}^{[1]}
   \end{array} \right).
\]

Then similar arguments as those used for (\ref{fk-est})
leads to the upper bound
\begin{align*}
    \left|u_k^{[j]}(\tau)\right|\le
    K\binom{k-1+1/m}{k}|1-\tau|^{k-2m-2},
    \qquad ( 1\le j\le 2(m-1), k\ge 0),
\end{align*}
for $\tau\in\mathscr{T}$, where the constant $K$ is easily
tuned according to the initial conditions.

\paragraph{Equality of two random median-of-($2t+1$) BSTs.} The
linear differential system of $4t+2$ equations of (\ref{de-m2tp1})
is
\begin{align*}
    U_j'(Z)&=U_{j+1}(Z),\quad 1\le j\le 2t, \\
    U_{2t+1}'(Z)&=(1-Z)^{-(2t+1)}U_{2t+2}(Z), \\
    U_j'(Z)&=U_{j+1}(Z),\quad 2t+2\le j\le 4t+1,  \\
    U_{4t+2}'(Z)&=\frac {(2t+1)!^2}{t!^4}\,\rho
    \sum_{0\le i_1, i_2\le t}
    \mu(i_1,i_2)(1-Z)^{2t-i_1-i_2} U_{2t+1-i_1}(Z)U_{2t+1-i_2}(Z),
\end{align*}
where
\[
    \mu(i_1,i_2):=\frac {(-1)^{i_1+i_2} \,t!^4}
    {i_1!i_2!(t-i_1)!^2(t-i_2)!^2}.
\]
Let $U_j(Z)=\sum_{k\ge 0} u_k^{[j]}(\tau)Z^{k-j-1}$ for $1\le j\le
4t+2$, where
\[
    u_0^{[j]}=(-1)^{j-1} j!\frac{(4t+3)!\,t!^4}{\rho\,(2t+1)!^4}
    \qquad(1\le j\le 2t+1).
\]
Then coefficient matrix $\mathbf{A}_k=k\mathbf{I}_{4t+2\times 4t+2}-\mathbf{M}$ in
(\ref{A-infinte-sys}), $k\ge 4t+2$, is given by
\[
    \mathbf{M}=\left(
    \begin{array}{ccccccccc}
             2 & 1 & 0 & \cdots &&&&&0 \\
             0 & 3 & 1 & \cdots & & &&&0 \\
             \vdots & \ddots & \ddots & \ddots &&&&&\vdots \\
             && 0 & 2t+1 & 1 &&&&\\
             &&& 0 & 2t+2 &1 &&&\\
             &&&& 0 & 2t+3 & 1 &&\vdots \\
             &&&&\vdots &\ddots &\ddots &\ddots &0 \\
             &&&&&& 0 & 4t+2 & 1 \\
             0 & 0 & \cdots & 0 &\frac {2(4t+3)!}{(2t+1)!}
             &&\cdots & 0 & 4t+3 \\
     \end{array}
     \right),
\]
and the vector-valued function $\mathbf{g}_k$ by
\[
    \mathbf{g}_k=\left(\begin{array}{c}
        0 \\ \vdots \\ 0\\
        \ds\sum_{0\le j<k} \binom{2t+k-j}{k-j}u_j^{[2t+2]}\\
        0\\ \vdots \\ 0 \\
        \ds\frac {(2t+1)!^2}{t!^4}\,
        \rho {\mathcal H}\left(u_j^{[t+1+\ell]}\biggr|
        \begin{array}{c}
            0\le j<k \\
            0\le\ell\le t
        \end{array}\right)
    \end{array} \right),
\]
where
\begin{align*}
    &{\mathcal H}\left(u_j^{[t+1+\ell]}\biggr|
    \begin{array}{c}
        0\le j<k \\
        0\le\ell\le t
    \end{array}\right)\\
    &\qquad=\mu(0,0)\sum_{\substack{0\le \ell \le k-j\\
    1\le j\le\min\{k,2t\}}}(-1)^j
    \binom{2t}{j}u_k^{[2t+1]}u_{k-j-\ell}^{[2t+1]}\\
    &\qquad\quad+\sum_{\substack{1\le s\le\min\{k,2t\}
    \\ 0\le i\le\min\{s,t\} }}\mu(i,s-i)\!\!
    \sum_{\substack{0\le \ell \le k-j\\ 1\le j\le\min\{k-s,2t-s\}}}
    \!\!(-1)^j\binom{2t-s}{j}u_k^{[2t+1-i]}u_{k-s-j-\ell}^{[2t+1+i-s]}.
\end{align*}

The the same method of proof used for (\ref{fk-est}) yields the
upper bound ($r_0=2$ or $6t+6$)
\begin{align*}
    \left|u_k^{[j]}(\tau)\right|
    \le C(1+k)^{-1/2}|1-\tau|^{k-r_0},
    \qquad ( 1\le j\le 4t+2,k\ge 0),
\end{align*}
uniformly for $\tau\in\mathscr{T}$, where the constants $C$ and $K$
are easily tuned according to the initial conditions.

\subsection*{A2.
Proof of the incompatibility of the resonance $r=2$
for random median-of-$(2t+1)$ BSTs}

Since the resonance $r=2$ does not depend on $t$, the
incompatibility of the resonance $r=2$ can be directly checked,
which we now do. Let $U(Z) := F(z)$, where $F$ satisfies the DE
(\ref{de-m2tp1}) and $Z=1-z/\rho$. Then the DE (\ref{de-m2tp1}) can
be rewritten as
\begin{align} \label{A-FringeDE}
    \left ( (1-Z)^{2t+1} U^{(2t+1)}(Z)\right)^{(2t+1)}
    =C_{t,\rho}\left(\left
    ( (1-Z)^{t}U^{(t)}(Z)\right)^{(t)}\right)^2,
\end{align}
where all derivatives are with respect to $Z$ and $C_{t,\rho} :=
(2t+1)!^2\rho/t!^4$.

Consider the formal Laurent expansion $f(Z)=\sum_{k\ge 0} u_k
Z^{k-\alpha}$. Then for any $s\in\mathbb{N}$, we have
\begin{align} \label{A-f-coeff}
    \left ( (1-Z)^{s}f^{(s)}(Z) \right)^{(s)}
    &= \sum_{k\ge 0} \falling{(k-\alpha-s)}{s}Z^{k-2s-\alpha}
    \sum_{0\le j\le s} (-1)^j\binom{s}{j}
    \falling{(k-\alpha-j)}{s}u_{k-j},
\end{align}
where $u_j:=0, j<0$. Substituting this into (\ref{A-FringeDE}), we
have
\begin{align*}
    & \sum_{k\ge 0}\falling{(k-\alpha-(2t+1))}
    {2t+1}Z^{k-4t-2\alpha}\sum_{0\le j\le 2t+1} (-1)^j\binom{2t+1}{j}
    \falling{(k-\alpha-j)}{2t+1}u_{k-j}\\
    &\qquad=C_{t,\rho}\sum_{k\ge 0}
    Z^{k-4t-2-\alpha}\sum_{0\le \ell\le k} \chi_k \chi_{k-\ell},
\end{align*}
where
\[
     \chi_k=\falling{(k-\alpha-t)}{t}
    \sum_{0\le j\le t} (-1)^j\binom{t}{j}
    \falling{(k-\alpha-j)}{t}u_{k-j}.
\]
Equating the dominant term (with $k=0$) leads to the obvious
solution $\alpha=2$. Consider now the relation
\begin{align*}
    & \falling{(k-\alpha-(2t+1))}{2t+1}
    \sum_{0\le j\le 2t+1} (-1)^j\binom{2t+1}{j}
    \falling{(k-\alpha-j)}{2t+1}u_{k-j}
    =C_{t,\rho}\sum_{0\le \ell\le k} \chi_k \chi_{k-\ell}.
\end{align*}
For $k=0$, we get $\rho u_0=(4t+3)!t!^4/(2t+1)!^4$, and
for $k=1$, we get $u_1=-2(t+1)^2u_0/(6t+5)$. Now
for $k=2$, we have
\begin{align*}
    0\cdot u_2 &=\left(-\frac {(4t+1)!}{(2t)!}
    \cdot\falling{0}{2t+1}+2(-1)^{t+1}\frac{(2t+1)!(2t-1)!}{(t-1)!}
    u_0\cdot\falling{0}{t}\right)u_2\\
    &= C_{t,\rho} (2t)!^2\left(
    \left((2t+1)(t+1)\binom{t}{2}+t^2(t+1)^2\right)u_0^2+u_1^2
    -t(4t+3)u_0u_1\right)\\
    & \mbox{}+(4t+1)!(2t+1)^2u_1-(4t+1)!(2t+1)(2t+2)
    \binom{2t+1}{2}u_0\\
    &=-\frac {(4t+2)!(t+1)}{4(6t+5)^2}u_0
    \left(216t^4+522t^3+437t^2+141t+12\right)\not=0,
\end{align*}
since $t\ge1$. This proves the incompatibility of the
resonace $r=2$ for all $t\ge1$.

\subsection*{A3. Asymptotics of the moment generating function}
We prove (\ref{MGF-X}), starting from Hankel's integral
representation of the Gamma function
\[
    \frac1{\Gamma(w)} = \frac1{2\pi i}\int_{\mathcal{H}_0}e^s s^{-w}
    \dd s \qquad(w\in\mathbb{C}),
\]
where $\mathcal{H}_0$ starts at $-\infty$, encircles the origin once
counter-clockwise and returns to its starting point. For
definiteness, we may take
\[
    \mathcal{H}_0 = \{s=xe^{\pm i\pi}\,:\, R_0\le x<\infty\}
    \cup \{s=R_0 e^{i\theta}\,:\, -\pi \le \theta\le \pi\}
    \qquad(R_0>0).
\]
This gives
\[
    M(z) = \frac1{2\pi i}\int_{\mathcal{H}_0} e^s s^{-1}
    A(z/s^{\alpha-1})\dd s,
\]
where $A(z)$ satisfies the DE (\ref{Az}). Note that $M$ is an
entire function of order $1/v>1$ and of type $\rho^{-1/v}$.

Let $z=|z|e^{i\varphi}$, $|z|>0$ and $|\varphi|<v\pi/2$, where
$v=(\sqrt{17}-3)/2$. The condition on $\arg z$ implies that the
dominant singularity $s=(z/\rho)^{1/v}$ of the integrand lies in the
half-plane $\Re(s)>0$ (in which $e^s\to\infty$ with $z$). On the
other hand, if $|\arg(-z)|<\pi- v\pi/2$, then one expects that
$M(z)\to0$ with $z$, but the exact determination of the rate is more
delicate. The situation here is similar to the Mittag-Leffler
function $\sum_{j\ge0} z^j/\Gamma(aj+1)$; see \cite[Ch.\
18.1]{EMOT55}.

The change of variables $z/s^v\mapsto s$ gives
\[
    M(z) = \frac1{2\pi iv}
    \int_{\mathcal{H}_1} e^{z^{1/v}
    s^{-1/v}}s^{-1}A(s) \dd s,
\]
where $\mathcal{H}_1$ is the cut circle described by
\[
    \mathcal{H}_1 = \{s=xe^{i\varphi\pm iv\pi}\,:\,0\le x\le
    R_1\}\cup
    \{s=R_1e^{i\varphi+iv\theta}\,:\,-\pi\le\theta\le\pi\}.
\]
Here $0<R_1<\rho$.  We then approach in a way similar to  the
singularity analysis (see \cite{FO90}) by deforming the contour
$\mathcal{H}_1$ into $\mathcal{H}_2$, where $\mathcal{H}_2$ is of
the same shape as $\mathcal{H}_1$ but with larger radius for the
circular part $|s|=R_2 = \rho+\ve$ and avoiding the cut from
$s=\rho$ to $\infty$ (in the style of \cite{FO90}). Symbolically,
\begin{align*}
    \mathcal{H}_2 &= \{s=xe^{i\varphi\pm iv\pi}\,:\,0\le x\le
    R_2\} \\ & \cup\,
    \{s=R_2e^{i\varphi+iv\theta}\,:\,-\pi\le\theta\le\pi
    \,\,{\rm and}\,\,
    |\theta-\varphi/v|\ge \ve_z\}
    \\ &  \cup\, \Gamma_\rho,
\end{align*}
where $\ve_z = |z|^{-1/v}$ and $\Gamma_\rho$ is any contour joining
the two points $R_2 e^{-i\ve_z}$ and $R_2 e^{i\ve_z}$ and lying
inside the cut region described by other parts of $\mathcal{H}_2$.

The remaining analysis is then easy because the main contribution to
$M(z)$ comes from $\Gamma_\rho$ on which we can apply the local
expansion (\ref{Az-ae}) of $A(z)$, the other parts being negligible
\[
    M(z) = \frac{1}{2\pi iv}
    \int_{\Gamma_\rho}
    e^{z^{1/v}s^{-1/v}}s^{-1}A(s) \dd s
    +O\left(e^{\Re(z/(\rho+\ve))^{1/v}}\right).
\]
By making first the change of variables $\rho(1-s)\mapsto s$, using
the expansion (\ref{Az-ae}), and then another change of variables
$(z/\rho)^{1/v} s/v\mapsto s$, we deduce that
\begin{align*}
    M(z) &= \frac{e^{(z/\rho)^{1/v}}}{2\pi i}\int_{\Gamma_0}
    e^s \Biggl(3\left(\frac{z}{\rho}\right)^{1/v}
    s^{-2}+\frac95 s^{-1}
    + \sum_{2\le j\le 7} \left(\bar{c}_j(s)+\tilde{c}_j(s)
    \log\frac z\rho\right) \left(\frac{z}{\rho}\right)^{-j/v} \\
    &\qquad +\frac{936}
    {21875}\left(\frac{z}{\rho}\right)^{-5/v}
    s^4\log s +O\left(|z|^{-6/v}|s|^5|\log s|\right)\Biggr)\dd s,
\end{align*}
where $\Gamma_0$ denotes the transformed contour of $\Gamma_\rho$
and the $c_j'$'s are polynomials of $s$ whose exact values matter
less. Extending the contour to infinity and then evaluating the
individual terms by Hankel's integral representation of the Gamma
function, we obtain
\[
    M(z) = e^{(z/\rho)^{1/v}}
    \left(3\left(\frac{z}{\rho}\right)^{1/v} +\frac95
    -\frac{22464}{21875}\left(\frac{z}{\rho}\right)^{-5/v}
    +O\left(|z|^{-6/v}\right)\right),
\]
where we also used the formula
\[
    \frac1{2\pi i}\int_{\mathcal{H}_0} e^s s^4\log s \dd s
    = -\frac{\dd{}}{\dd x}
    \frac1{\Gamma(x)}\Biggr|_{x=-4} = -24.
\]
This completes the proof of (\ref{MGF-X}).

\end{document}